\documentclass{article}

\usepackage{arxiv}

\usepackage[utf8]{inputenc} 
\usepackage[T1]{fontenc}    
\usepackage{hyperref}       
\usepackage{url}            
\usepackage{booktabs}       
\usepackage{amsfonts}       
\usepackage{nicefrac}       
\usepackage{microtype}      
\usepackage{lipsum}		
\usepackage{graphicx}
\usepackage{natbib}
\usepackage{doi}
\usepackage{bbding}
\usepackage{array}
 \usepackage{array,multirow,graphicx}

\usepackage{amsmath}
\usepackage{diagbox}
\usepackage{color}
\usepackage{bm}
\usepackage{array}
\usepackage[utf8]{inputenc}
\usepackage{graphicx}
\usepackage{float}
\usepackage{amsmath}
\usepackage{caption}
\usepackage{subcaption}
\usepackage{multirow}

\usepackage{color}
\usepackage{booktabs}
\usepackage[ruled,vlined]{algorithm2e}
\usepackage{makecell}
\newcolumntype{I}{!{\vrule width 1.5pt}}
\usepackage{booktabs}

\usepackage{graphicx}  
\usepackage{soul}
\usepackage{subcaption}
\usepackage{mwe}
\usepackage{xcolor}
\usepackage{amssymb}

\usepackage{tikz}
\definecolor{fc}{HTML}{1E90FF}
\tikzset{fc/.style={black,draw=black,fill=fc,rectangle,minimum height=1cm}}
\definecolor{h}{HTML}{228B22}
\definecolor{bias}{HTML}{87CEFA}
\tikzset{h/.style={black,draw=black,fill=h,rectangle,minimum height=1cm}}
\tikzset{bias/.style={black,draw=black,fill=bias,rectangle,minimum height=1cm}}
\definecolor{anti-flashwhite}{rgb}{0.95, 0.95, 0.96}
\definecolor{almond}{rgb}{0.98, 0.91, 0.71}

\title{Physics-informed Discretization-independent \\ Deep Compositional Operator Network}

\author{
  Weiheng Zhong \\
  Department of Civil and Environmental Engineering\\
  University of Illinois at Urbana-Champaign\\
  Champaign, Illinois\\
  \texttt{weiheng4@illinois.edu} \\
   \And
 Hadi Meidani \\
  Department of Civil and Environmental Engineering\\
  University of Illinois at Urbana-Champaign\\
  Champaign, Illinois \\
  \texttt{meidani@illinois.edu} \\
}
 
\begin{document}
\maketitle

\begin{abstract}
Solving parametric Partial Differential Equations (PDEs) for a broad range of parameters is a critical challenge in scientific computing. To this end, neural operators, which \textcolor{black}{predicts the PDE solution with variable PDE parameter inputs}, have been successfully used. However, the training of neural operators typically demands large training datasets, the acquisition of which can be prohibitively expensive. To address this challenge, physics-informed training can offer a  cost-effective strategy. However, current physics-informed neural operators face limitations, either in handling irregular domain shapes or in \textcolor{black}{in generalizing to various discrete representations of PDE parameters}. In this research, we introduce a novel physics-informed model architecture  which can generalize to \textcolor{black}{various discrete representations of PDE parameters} and  irregular domain shapes. Particularly, inspired by deep operator neural networks, our model involves a discretization-independent learning of parameter embedding repeatedly, and this  parameter embedding is integrated with the response embeddings through multiple compositional layers, for more expressivity.  Numerical results  demonstrate the accuracy and efficiency of the proposed method. All the  codes and data related to this work are available on GitHub: \textcolor{blue}{\href{https://github.com/WeihengZ/PI-DCON}{https://github.com/WeihengZ/PI-DCON}} .
\end{abstract}


\section{Introduction \label{Sec.Intro}}

Partial Differential Equations (PDEs) are used to describe the behaviors of systems in various fields such as physics \cite{phy_app1,phy_app2,phy_app3,phy_app4,phy_app5}, chemistry, and biology \cite{pde_app}. The predominant numerical approach to solve  PDEs to calculate the system response is the Finite Element Method (FEM) \cite{fem}, which discretizes the continuous domain in which the response is calculated, into a mesh. For complex problems with significant  nonlinearities calculating the FEM solution is computationally expensive, especially in design optimization or uncertainty quantification tasks requiring  repetitive  solution of parametric PDEs for a  range of conditions and parameters \cite{TO}.

Machine learning techniques have been introduced to accelerate the process of solving parametric PDEs by learning a neural operator, which serves as a mapping from the PDE parameters to the system response (solution of the PDE) \cite{DON}. If neural operators are successfully trained on a dataset, they  can generalize to new, unseen PDE parameters. This means that for a new PDE parameter, the solution can be calculated  via a single forward pass on the trained neural network, with minimal computational time. 

Recently,  various neural operator architectures have been proposed for operator learning in a data-driven (supervised learning) approach \cite{DON, FNO, NO, GNO}. As the pioneering work in this area, Deep Operator Networks (DeepONet) \cite{DON} were developed, based on an inspiration from the universal approximation theorem for operators. DeepONet approximates nonlinear operators by effectively learning a collection of basis functions and coefficients, in the form of neural networks. However, the function input to DeepONet, i.e. the PDE parameters, are structured in the form of a finite-dimensional vector with a fixed dimension. This makes DeepONet to be  a \textcolor{black}{ discretization-specific} method \cite{NO}. \textcolor{black}{For instance, consider employing sensors on $m$ distinct locations to observe a PDE parameter, leading to  an $m$-dimensional representation of the PDE parameter. Then, a DeepONet is trained and only usable for an exactly $m$-dimensional parameter input,  and won't be applicable for cases where more or fewer sensors are used, or when  discretizations with a different size are used. This restriction limits the utility of neural operators for cases where various discretizations may be used.}

\textcolor{black}{However, in most engineering applications, we may use different mesh resolution in different location of the domain, focusing on different aspect of the system performance evaluation.} To overcome the limitations of DeepONet, numerous \textcolor{black}{"discretization-independent"} neural operators have been introduced. \textcolor{black}{These discretization-independent methods can predict the PDE solutions with various discrete representation of the PDE parameters on different meshes.} The Graph Neural Operator (GNO) \cite{GNO} adopts graph neural networks (GNNs) architecture for operator learning, i.e. by considering graph structures for model's input and output. The Fourier Neural Operator (FNO) \cite{FNO} leverages the Fourier transform to learn mappings in the spectral domain, capturing global dependencies and showing superior performance in rectangular domain shapes. Geo-FNO \cite{geo_FNO} extends FNO to irregular meshes by learning a mapping from an irregular mesh to a uniform mesh. Inspired by FNO, Wavelet Neural Operator (WNO) \cite{WNO} uses the Wavelet Transformation, instead of the Fourier transformation, to better model  signals with discontinuity and spikes, and also handle non-square domains. Also, the Low-Rank Neural Operator (LRNO) \cite{NO} boosts the efficiency and scalability of neural operator models, using a low-rank representation of the kernels in neural operators. 

While these "discretization-independent" methods have shown promise in various applications, they all are data-intensive \cite{PI-deepOnet}. Training neural operators in a data-driven (supervised) way necessitates collection of a large number of parameter-solution pairs as the training set, which is obtained using  repeated  costly high-fidelity (FEM) simulations. This is prohibitively expensive in complex industrial applications. Under this circumstance,  physics-informed neural operators emerge as an effective solution for constructing neural operators. These methods draw inspiration from physics-informed neural networks (PINNs) \cite{pinn}, and integrate the governing PDEs directly into  the training process. This allows for a data-free training  without need for  FEM-based training data. Among the first physics-based methods, parametric physics-informed neural networks \cite{PINN_bc, PINN_param, PINN_geo} considered  various elements of a PDE problem, such as PDE parameters \cite{PINN_param, PNN_param_gen}, boundary conditions \cite{PINN_bc}, and even the geometry of the domain \cite{PINN_geo}, as additional  feature vectors.  These feature vectors provide additional inputs to the model, allowing it to learn the dependence of the PDE solution on PDE parameters. However, this approach is limited to cases where an analytical  representation of varying conditions, such as boundary conditions and geometries, are available, and cannot be effective where these conditions are available in discrete form with many nodes.

Another class of data-free methods are neural operators that are trained using a physics-informed loss function \cite{PPINN1,PPINN2,PPINN3, PPINN4}. Among them, physics-informed DeepONet (PI-DeepONet) \cite{PI-deepOnet} extends the original DeepONet framework by embedding physical laws directly into the loss function which can be effectively calculated in any irregular domain shape. The incorporation of physical laws also enhances the learning process but, even with this advancement, physics-informed DeepONet cannot  generalize to   parameter discretizations with different size. On the other hand, the physics-informed Fourier Neural Operator (PI-FNO) \cite{PI-FNO} builds upon the original FNO framework, and offers generalization across different  parameter dicretizations while also incorporating physics-informed training methods. However, the reliance of PI-FNO on Fast Fourier Transform (FFT) restricts its application to rectangular domains with uniform meshes as FFT works  on uniform sampling of the signal in the spatial domain \cite{geo_FNO}. Another notable development is the physics-informed Wavelet Neural Operator (PI-WNO) \cite{PI-WNO}, which implements physics-informed training for WNO \cite{WNO} using stochastic projection. Although PI-WNO can handle various PDE parameter discretizations, the hidden kernels are learned based on the architecture of convolution neural networks \cite{cnn}, which prevents the model from learning in complex domain geometries. 

In our study, we \textcolor{black}{seek to overcome the aforementioned limitations related to } data requirement, discretization-dependence and irregular geometries  by introducing an innovative neural operator architecture. Our Physics-informed Deep Compositional Operator Network (PI-DCON) model is capable of generalizing across different \textcolor{black}{ discrete representations of  PDE parameters}, including those in irregular domain shapes. In this study, the irregular domain shapes are defined as those that does not belong to regular geometric forms such as rectangles, circles, or polygons. The differences between our proposed model versus other existing works are summarized in Table \ref{table.archi_compare}. To the best of our knowledge, our paper is the first attempt to develop a \textcolor{black}{discretization-independent} neural operator that can handle complex geometries without need for any training data.  

\begin{table*}[ht]
\begin{center}
\caption{Differences of Physics-informed Deep Compositional Operator Network versus other existing models.}
\begin{tabular}{c c c c}
\hline
{} & Data-free & \makecell{Generalize across PDE \\ parameter representations} & \makecell{Handle irregular \\ domain shapes} \\
DeepONet & \XSolidBrush & \XSolidBrush & \CheckmarkBold \\
FNO, WNO ,LRNO & \XSolidBrush & \CheckmarkBold & \XSolidBrush  \\
GNO, Geo-FNO & \XSolidBrush & \CheckmarkBold & \CheckmarkBold  \\
PI-DeepONet & \CheckmarkBold & \XSolidBrush & \CheckmarkBold \\
PI-FNO , PI-WNO  & \CheckmarkBold & \CheckmarkBold & \XSolidBrush \\
PI-DCON & \CheckmarkBold & \CheckmarkBold & \CheckmarkBold \\
\hline
\label{table.archi_compare}
\end{tabular}
\end{center}
\end{table*}

The remainder of this paper is organized as follows. In Section~\ref{sec.bg}, we  briefly introduce the problem settings and the technical backgrounds of PINNs and PI-DeepONet. Section~\ref{sec.methodology} briefly introduced our model architectures and the high-level ideas behind our model. Finally, a detailed analysis of the performance evaluation of the proposed methods and conclusions are included in Sections~\ref{sec.results} and \ref{sec.conclusions}.

\section{Technical background \label{sec.bg}}

\subsection{Problem Setting \label{Subsec.PS}}

In this study, our goal is to develop an efficient machine-learning based solver for solving parametric PDEs \cite{pde_app} which are formulated by: 
\begin{equation}
    \begin{aligned}
    \mathcal{N}_{\bm{x}} [u(\bm{x}), k(\bm{x})] &= 0, \quad & \bm{x}  \in \Omega, \\
    \mathcal{B}_{\bm{x}}[u(\bm{x})] &= g(\bm{x}),  \quad  & \bm{x} \in \partial \Omega,
\end{aligned}
\label{eq.pde}
\end{equation}
where $\Omega$ is a $d$-dimensional physical domain in $R^d$, $\bm x$ is a $d$-dimensional spatial coordinate, $\mathcal{N}_{\bm x}$ is a general differential operator, and $\mathcal{B}_{\bm{x}}$ is a boundary condition operator acting on the domain boundary $\partial \Omega$.  Also, $k(\bm{x})$ refers to the parameters of the PDE, which can include the  coefficients and forcing terms in the governing equation; $g(\bm{x})$  denotes the  boundary conditions;  and  $u(\bm{x})$ is the solution of the PDE at the given parameters and boundary conditions. We will use neural network models to approximate the solution operator $\mathcal{M}: \{k(\bm{x}),g(\bm{x})\} 
\rightarrow u(\bm{x})$. 

\subsection{Physics-informed Neural Networks \label{subsec.pinn}}

Given actual parameters and boundary conditions $\{k^*(\bm{x}),g^*(\bm{x})\}$ and the solution operator $\mathcal{M}$ defined by equations \ref{eq.pde}, we denote $u^{*}(\bm{x}) = \mathcal{M} (k^{*}(\bm{x}),g^{*}(\bm{x}))$ as the unique ground truth. The solution $u^{*}(\bm{x})$ can be approximated by a neural network $u_{\theta}(\bm{x})$ using physics-informed training \cite{pinn}, where $\theta$ denotes the neural network parameters, i.e. weights and biases. By incorporating the physics laws into training loss, the total training loss of PINNs $\mathcal{L}_{\text{PINN}}$ is formulated as:
\begin{equation}
    \mathcal{L}_{\text{PINN}} = \mathcal{L}_{\text{PDE}}  + \alpha \mathcal{L}_{\text{BC}},
\label{eq.pinn_loss}
\end{equation}
where
\begin{equation}
    \begin{aligned}
    \mathcal{L}_{\text{PDE}} &= (\mathcal{N}_{\bm{x}} [u_{\theta}(\bm{x}), k^*(\bm{x})])^2,  &\bm{x} \in \Omega, \\
    \mathcal{L}_{\text{BC}} &=  (\mathcal{B}_{\bm{x}}[u_{\theta}(\bm{x})] - g^*(\bm{x}))^2, &\bm{x} \in \partial \Omega,
        \end{aligned}
\label{eq.pinn_pde_loss}
\end{equation}
where $\alpha$ is the trade-off coefficient between the PDE residual loss term and the boundary conditions loss function. The optimal neural network parameters $\theta$ are found by minimizing the total training loss $\mathcal{L}_{\text{PINN}}$ with exact derivatives computed using automatic differentiation \cite{AutoD}. However, each well-trained neural network can only approximate a single solution with respect to one realization of the parameters, and will not be suitable for solving parametric PDEs with varying realizations of the parameters. 

\subsection{Learning operators with physics-informed DeepONet \label{PI-DeepONet}}

An alternative approach to solve parametric PDEs is approximating the operator $\mathcal{M}$ directly with a neural network $u_{\theta}(\bm{x}; k(\bm{x}'), g(\bm{x}'))$, where $\bm{x}$ represents the locations at which the solution $u$ is calculated, and $\bm{x}'$ denote the locations on which the values of parameters are available. In order to set up the training loss, let us consider $M$ to be the total number of  realizations of parameter functions $k(\cdot)$ and $g(\cdot)$, each obtained at locations $\bm{x}'$. Then the physics-informed training loss of neural operator $\mathcal{L}_{\text{PINO}}$ is given by 
\begin{equation}
    \mathcal{L}_{\text{PINO}} = \mathcal{L}_{\text{PDE}}  + \alpha \mathcal{L}_{\text{BC}},
\label{eq.pino_loss1}
\end{equation}
where
\begin{equation}
    \begin{aligned}
    \mathcal{L}_{\text{PDE}} &= \frac{1}{N} \sum_{i=1}^N (\mathcal{N}_{\bm{x}} [u_{\theta}(\bm{x}; k_i(\bm{x}'), g_i(\bm{x}')), k_i(\bm{x}')])^2,  &&\bm{x} \in \Omega, \\
    \mathcal{L}_{\text{BC}} &= \frac{1}{N} \sum_{i=1}^N (\mathcal{B}_{\bm{x}}[u_{\theta}(\bm{x}; k_i(\bm{x}'), g_i(\bm{x}'))] - g_i(\bm{x}'))^2, &&\bm{x} \in \partial \Omega,
\end{aligned}
\label{eq.pino_loss}
\end{equation}
where $N$ is the number of PDE parameter realizations. For any new realization of the parameters, the well-trained neural operator can predict the corresponding solution directly. A popular choice of model architecture to approximate the PDE operator is DeepONet \cite{DON}. DeepONet is composed of two separate neural networks referred to as the "branch net" and "trunk net", respectively. Both the branch net and trunk net are simply multilayer perceptrons.

Without  loss of generality, let us consider a DeepONet that approximates the PDE operator mapping from boundary conditions to solutions $\mathcal{M}: g(\bm{x}') \rightarrow u(\bm{x})$. As shown in Figure~\ref{fig.inspirationa}, the solution is calculated using  the \emph{parameters embedding} and the \emph{response embedding}. Specifically, the input to the branch net is $[g(\bm{x}'_1), g(\bm{x}'_2), ..., g(\bm{x}'_m)]$  which is  function $g(\bm{x}')$ evaluated at a collection of fixed locations $\{\bm{x}'_i\}_{i=1}^m$. The output of the branch net is a $q$-dimensional parameters embedding $\bm{b}=[b_1,b_2,..., b_q]$. The trunk net takes the continuous coordinates $\bm{x}$ as input, and outputs a $q$-dimensional coordinate embedding $\bm{t} = [t_1,t_2,..., t_q]$. The outputs of the branch net and trunk net are merged together by dot product to produce  the solution $u$  at location $\bm{x}$ by the following equation:
\begin{equation}
    u_{\theta}(\bm{x}, g(\bm{x}')) = \sum_{i=1}^q b_i t_i(\bm x).
    \label{eq.DeepONet}
\end{equation}

The parameters of the DeepONet $\theta$ can be optimized by the training loss in equation \ref{eq.pino_loss}. 

\section{Methodology}\label{sec.methodology}

\begin{figure}[!ht]
     \centering
     \begin{subfigure}[b]{0.43\textwidth}
         \centering
         \includegraphics[height=2.8cm]{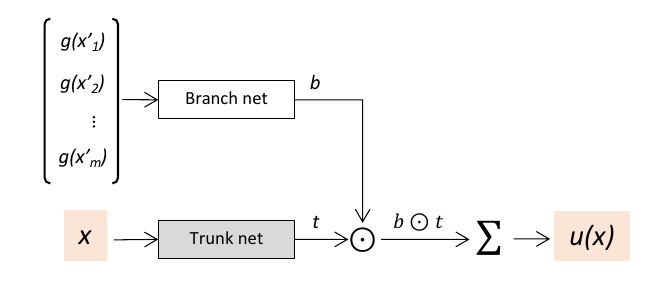}
         \caption{DeepONet}
         \label{fig.inspirationa}
     \end{subfigure}
     \begin{subfigure}[b]{0.55\textwidth}
         \centering
         \includegraphics[height=2.8cm]{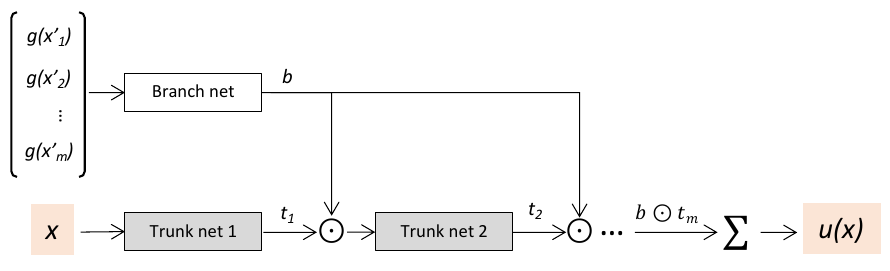}
         \caption{Deep Compositional Operator Network}
         \label{fig.inspirationb}
     \end{subfigure}
        \caption{\footnotesize The connection between our model architecture and DeepONet is shown. The architecture on the left is DeepONet and the one on the right is our proposed model. Comparing two model architectures, we can observe that our model is a compositional version of the DeepONet.}
\end{figure}

Our proposed model architecture can be considered as a modified version of DeepONet. Let us first begin from the Universal Approximation Theorem for Operator learning \cite{DON}, which states that any functional operator can be approximated using  parameter  and coordinate embeddings. Specifically, as shown in Figure~\ref{fig.inspirationa}, in DeepONet, the operator is approximated as the inner product $\bm b \cdot \bm t $, or when written differently, as 
\begin{equation}
    \mathcal{M} \approx \bm b \cdot \bm t = \text{sum} \{\bm{b} \odot \bm t(\bm{x})\},
    \label{eq.theorem_DeepONet}
\end{equation}
where $\odot$ refers to component-wise multiplication, and $\text{sum}(\cdot)$ returns the summation of the components of a vector. 

Inspired by this formulation, we conjecture that a compositional framework to combine the parameter and response embeddings  can lead to a more expressive model that can approximate more complex functional operators. In particular, we introduce the Deep Compositional Operator Network (DCON) which approximates the PDE operator as follows
\begin{equation}
    \mathcal{M} \approx \text{sum} \{\bm b \odot ... \ N(\bm b \odot N(\bm b \odot \bm t(\bm{x})))\},
    \label{eq.theorem_DCON}
\end{equation}
where $N: \mathbb{R}^q \rightarrow \mathbb{R}^q$ is a learnable nonlinear mapping (see Figure \ref{fig.inspirationb}). 

Without loss of generality, let us consider a PDE operator mapping from boundary conditions to the PDE solution, denoted by $\mathcal{M}: g(\bm{x}') \rightarrow u(\bm{x})$. In DCON, similarly to DeepONet, the boundary conditions is encoded using a finite number of sampled boundary conditions $[\bm{x}'_1, \bm{x}'_2, ..., \bm{x}'_m]$, and the corresponding function values $[g(\bm{x}'_1), g(\bm{x}'_2), ..., g(\bm{x}'_m)]$.  The  architecture of DCON consists of  a "branch network" and a "solution network". The branch network  maps the given boundary conditions on each of the available $m$ boundary locations to a higher dimensional space, $\mathcal{S}_\text{param}$, using a multi-layer perceptron. Then, a Max-pooling layer is applied to this high-dimensional representation to capture, in each dimension of $\mathcal{S}_\text{param}$, the most important feature across the $m$ locations on the boundary. The final output of the branch network is a global parameter embedding $\bm b$ in the high-dimensional space $\mathcal{S}_\text{param}$, whose dimension is independent of the number of boundary  points, $m$. This setup relaxes the requirement for all training and test data to include the same number of boundary locations, $m$. \textcolor{black}{In contrast to FNO, which aggregates spatial information after mapping spatial features to the frequency domain, we employ a simple yet effective pooling operation to directly aggregate parameter information in the original spatial domain. By avoiding the use of FFT, our model offers an efficient framework for physics-informed training using Auto Differentiation \cite{AutoD}.}

The inputs to the solution network are the coordinates of a  collocation point, denoted by $\bm{x}$, and the parameter embedding, $\bm b$, based on which the solution  at the collocation points $u(\bm{x})$ is calculated. As shown in equation \ref{eq.theorem_DCON}, the functional mapping $\mathcal{M}$ is approximated with a number of operator layers $H$ formulated as following:
\begin{equation}
    v^{l+1} = H(v^l) := \bm b \odot N(v^l),
\end{equation}
where $v^l$ represents the output of $l$-th operator layer. In our proposed model, we simply use a linear layer with a nonlinear activation function as the mapping $N$. Hence, the PDE operator can be approximated by our proposed model using the following equation:
\begin{equation}
    \mathcal{M} \approx \text{sum} \{\bm b \odot  ... (\bm b \odot (W^2_L \ \sigma(\bm b \odot (W^1_L \ \sigma (W_{\bm t}  \bm{x} + B_{\bm t}) + B^1_L)) + B^2_L)  )\},
    \label{eq.formula_model_arch}
\end{equation}
where  $\sigma$ represents the nonlinear activation function, $W_{\bm t} \in R^{q \times d}, B_{\bm t} \in R^{q}$ are  trainable parameters of the mapping $\bm t$, and $W^j_L \in R^{q \times q}, B^j_L \in R^{q}$ are  trainable parameters of the $j$-th mapping $N$. These parameters of the model are  optimized by minimizing the physics-informed training loss as shown in equation \ref{eq.pinn_loss}. 

\textcolor{black}{Based on the Universal Approximation Theorem, a single operator layer can approximate any continuous functional operator. Therefore, using an architecture of stacked operator layers, we can achieve at least the same approximation capability as DeepONet. Regarding model complexity, our model architecture is also comparable to DeepONet. For a hidden feature dimension of $q$, DeepONet with $L$ hidden layers has model complexity of $\mathcal{O}(Lq^2 + q)$, considering the matrix-vector multiplication in each hidden layer and the element-wise multiplication in the last layer. In contrast, our model has similar model complexity of $\mathcal{O}(Lq^2 + Lq)$ due to the additional element-wise multiplication in each hidden layer.}

\begin{figure}[H]
    \centering
    \includegraphics[width=16cm]{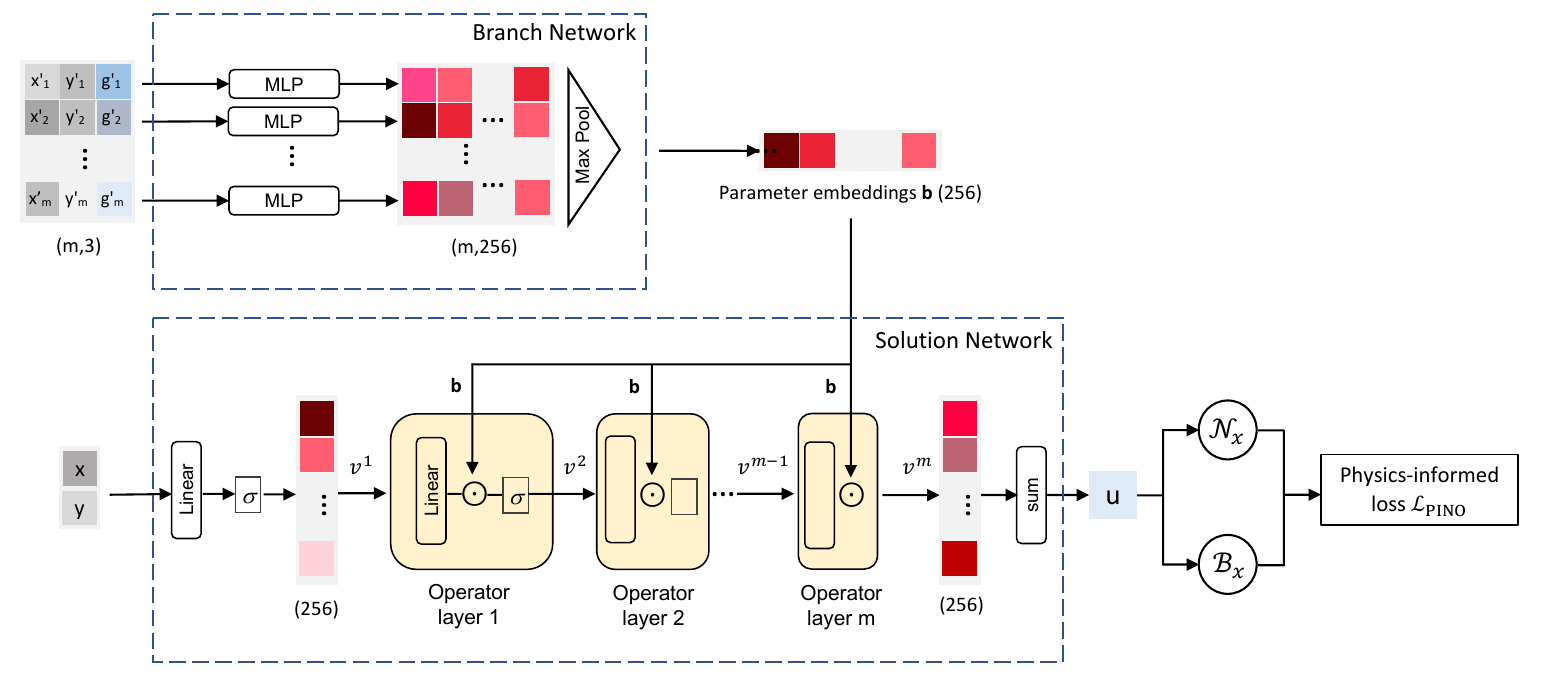}
    \caption{\footnotesize The architecture of the physics-informed Deep Compositional Operator Network is shown. The grey blocks represent the coordinates and the blue blocks represent the function values. The red blocks represent the hidden embeddings. The value of $b'_i$ is the boundary condition value evaluated on $(x'_i, y'_i)$. The value of $u$ is the solution value evaluated on $(x, y)$. It should be noted that no activation function is included in the last operator layer.}
    \label{fig.PIK}
\end{figure}

More details about the proposed model architecture are shown in Figure \ref{fig.PIK}. It should be noted that the proposed DCON architecture can also be used in a data-driven approach, where its parameters are estimated in a supervised learning process, as will be discussed in Section~\ref{sec.datadriven}.

\section{Numerical results} \label{sec.results}

In this section, we numerically evaluate the accuracy of the proposed model for solving parametric differential equations. We compare our models with data-driven DeepONet \cite{DON}, physics-informed DeepONet  \cite{PI-deepOnet}, and \textcolor{black}{an improved architecture of DeepONet reported in \cite{IDON}.} We use these default settings unless mentioned otherwise; We use the hyperbolic tangent function (Tanh) \cite{tanh} as our activation function to ensure the smoothness in high-order derivatives. For the architecture of DeepONet, we use three hidden layers of width 512. For our model architecture, we use three operator layers of width 512. \textcolor{black}{
Comparing different learning methods is challenging because of the numerous hyper-parameters involved. To address this, we use the grid search method \cite{grid-search} to select the appropriate hyper-parameters for each approach, ensuring the hyperparameters are properly set for different learning methods.} Adam is the default optimizer with the following default hyper-parameters: $\beta_1$ = 0.9 and $\beta_2$ = 0.999 \cite{adam}. Each model is trained on 70\% of the sampled PDE parameters and validated on 10\% of the sampled PDE parameters. Well-trained models predict the PDE solution for the remaining 20\% of the sampled PDE parameters. Two hyper-parameters used in model training: (1) the learning rate, with the possible values of learning rate are 0.001, 0.0005, 0.0002, and 0.0001; and (2) the coordinate sampling size ratio  in each epoch. This ratio determines the number of collocation points in each epoch as a function of the number of sampled PDE parameters. The possible values for this ratio are 0.3, 0.2, 0.1, and 0.05. Model training is performed on an NVIDIA P100 GPU using a batch size of 20. We collect the parameters of PDEs by Monte Carlo Simulation \cite{mcs} of stochastic processes \cite{sp} and derive the solution of PDEs based on the  Finite Element Method using Matlab \cite{matlab}. \textcolor{black}{We use the $L_2$ relative error ${|u_{\text{pred}}-u_{\text{fem}}|} / {|u_{\text{gt}}|}$ as our error measure, where $u_{\text{pred}}$ refers to the neural network prediction of PDE solution, and $u_{\text{fem}}$ is the PDE solutions computed by FEM.}

\subsection{Experiment Setups}

\subsubsection{A Darcy flow problem}

\begin{figure}[ht]
    \centering
    \includegraphics[width=\textwidth]{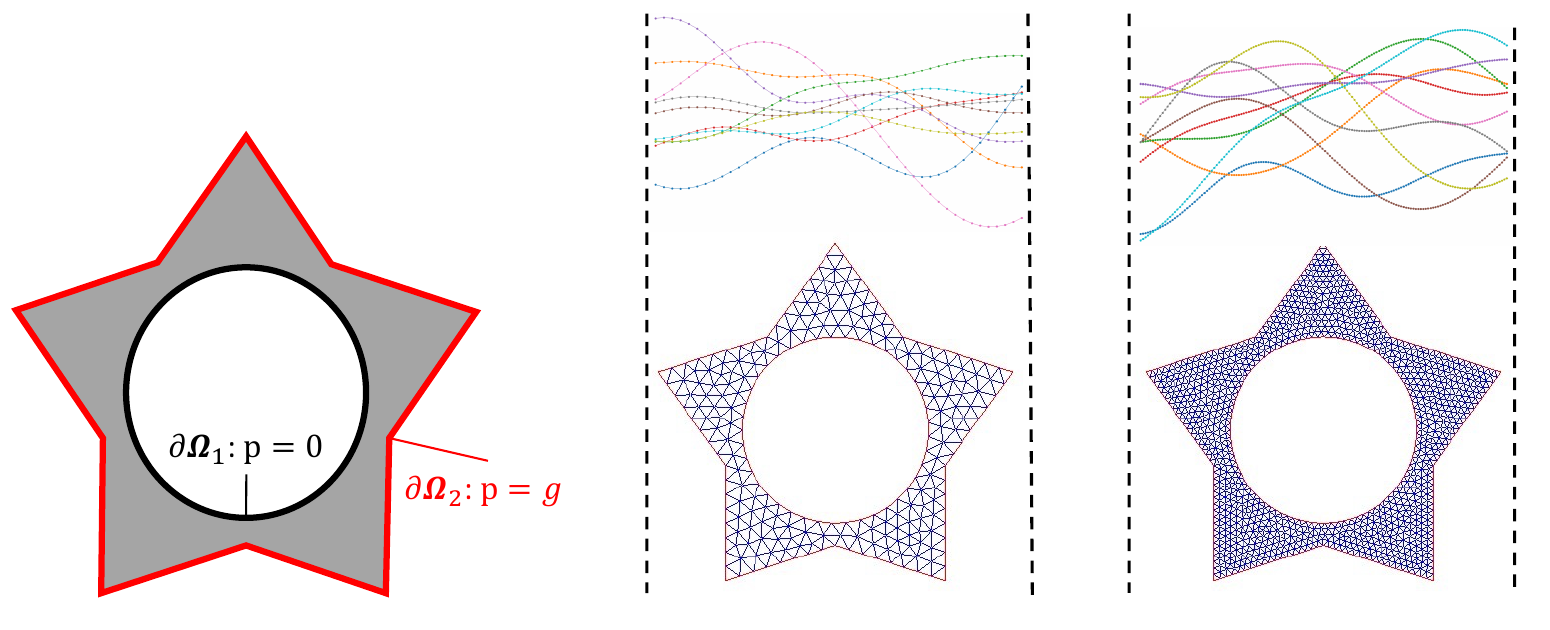}
    \caption{\footnotesize The boundary conditions of our Darcy flow experiment setting, and the coarsest and finest meshes used as discretizations to solve the PDE are shown. The samples of Gaussian processes as the boundary conditions $g(\bm{x})$ are also shown above the meshes. The dots on the curve are considered as the representation of the boundary condition values.}
    \label{fig.darcy_BC}
\end{figure}

Let us consider a two-dimensional Darcy flow problem in a pentagram-shaped domain with a hole inside, which is a benchmark problem studied by a data-driven neural operator used in \cite{NO_compare}. The steady state solution of the system is described by the following equation:
\begin{equation}
    \begin{aligned}
    - \nabla (k(x,y) \nabla p(x,y)) &= f(x,y), \quad (x,y) \in \Omega,
\end{aligned}
\end{equation}
where $p(x,y)$ is the pressure, $k(x,y)$ is the permeability field,  $f(x,y)$ is the source term, and $g(x,y)$ is the prescribed pressure on the domain boundaries. Using the same setting of \cite{NO_compare}, we consider $k(x,y)=1$ and $f(x,y)=10$. The boundary consists of two parts, i.e. $\partial \Omega = \partial \Omega_1 \cup \partial \Omega_2 $. First,  on the perimeter of the hole, denoted by $\partial \Omega_1$,  the pressure is set to be zero. Second,  on the outer edges of the pentagram, denoted by $\partial \Omega_2$, we impose a nonzero boundary condition function $g$. Hence, the boundary conditions of the problem is formulated as:
\begin{equation}
\begin{aligned}
        p(x,y) &= 0,  &(x,y) \in \partial \Omega_1, \\
        p(x,y) &= g(x,y), &(x,y) \in \partial \Omega_2.
\end{aligned}
\end{equation}
In this example, we consider the pressure on $\partial \Omega_2$ to follow a zero mean Gaussian process with a covariance kernel that is only dependent on the horizontal distance, i.e., 
\begin{equation}
    \begin{aligned}
    g(x,y) &\sim \mathcal{GP}(0, K(x,x')),\\
    K(x,x') &= \exp\left[-\frac{(x-x')^2}{2l^2}\right].
\end{aligned}
\label{eq.gp}
\end{equation}
 We employ this Gaussian process model of Eq.~\ref{eq.gp} with $l=1$ to generate $N=1000$ different imposed pressure functions on $\partial \Omega_2$. Then given these sampled functions, $\{g_i\}_{i=1}^N$ we solve the PDE  using Finite Element Method (FEM) \cite{pde_toolbox}. In order to demonstrate the \textcolor{black}{resolution-independence} feature of our model, for different sampled boundary functions,  we consider different mesh sizes in the domain. Specifically, for the $i$-th realization, with the boundary condition function $g_i$, let $M_{1,i}$ and $M_{2,i}$ be the number of sampled locations on $\partial \Omega_1$ and $\partial \Omega_2$, respectively, and $\bm{X}_{\text{BC}_1,i}$ and $\bm{X}_{\text{BC}_2,i}$ are the corresponding sets of coordinates at those sampled locations, i.e.
\begin{equation}
\begin{aligned}
    \bm{X}_{\text{BC}_1,i} = \{(x_j^{1,i}, y_j^{1,i})\}_{j=1}^{M^{1,i}}, \\
    \bm{X}_{\text{BC}_2,i} = \{(x_j^{2,i}, y_j^{2,i})\}_{j=1}^{M^{2,i}},
\end{aligned}
\end{equation}
Specifically, we choose $M_{1,i}$ and $M^{2,i}$ to be random numbers drawn from $[30,80]$ and   from $[100,300]$, respectively. Therefore, all the boundary information for a single data is collected in  $\bm{G}^i$  as:
\begin{equation}
    \bm{G}_i = \{\bm{X}_{\text{BC}_2,i}, g_i(\bm{X}_{\text{BC}_2,i})\},
\end{equation}
where $g_i(\bm{X}_{\text{BC}_2,i})$ is the boundary values evaluated at the coordinates $\bm{X}_{\text{BC}_2,i}$. Figure ~\ref{fig.darcy_BC} shows visualizations of a few sampled  boundary conditions together on the finest and coarsest meshes used in this study.

Given this setup, we seek to learn the nonlinear mapping $\mathcal{M}_{\text{Darcy}}$ that transforms a given boundary conditions $g(x,y)$ on $\partial \Omega_2$ to the pressure field $p(x,y)$ in the entire domain, i.e.,
\begin{equation}
    \mathcal{M}_{\text{Darcy}}: g(x,y) \rightarrow p(x,y).
\end{equation}
Specifically, for a 2D problem, the neural operator $p_{\theta}(x,y, \bm G)$  
is trained using the loss function
\begin{equation}
    \mathcal{L}(\theta) = \mathcal{L}_{\text{PDE}}(\theta) + \alpha_1 \mathcal{L}_{\text{BC}_1}(\theta) + \alpha_2 \mathcal{L}_{\text{BC}_2}(\theta),
    \label{eq.darcy_loss}
\end{equation}
where $\alpha_1$, $\alpha_2$ are the weighing hyperparameters, set to be $\alpha_1=\alpha_2=500$ to ensure similar orders of magnitude. $\mathcal{L}_{\text{PDE}}(\theta)$ and $\mathcal{L}_{\text{BC}}(\theta)$ are the PDE residual and the BC loss function, given by
\begin{equation}
\begin{aligned}
    \mathcal{L}_{\text{PDE}}  (\theta) &= \frac{1}{N} \sum_{i=1}^N \left[\frac{\partial^2 p_{\theta}(x,y, \bm G_i)}{\partial x^2} + \frac{\partial^2 p_{\theta}(x,y, \bm G_i)}{\partial y^2} + 10 \right] ^2,  &(x, y) \in \Omega, \\
    \mathcal{L}_{\text{BC}_1}(\theta) &= \frac{1}{N} \sum_{i=1}^N \left\{\frac{1}{M_{1,i}} \sum_{j=1}^{M_{1,i}} \left[ p_{\theta}(x^i_j,y^i_j, \bm G_i) - 0 \right]^2 \right\},  &(x^i_j, y^i_j) \in \bm{X}_{\text{BC}_1,i}, \\
    \mathcal{L}_{\text{BC}_2}(\theta) &= \frac{1}{N} \sum_{i=1}^N \left\{\frac{1}{M_{2,i}} \sum_{j=1}^{M_{2,i}} \left[ p_{\theta}(x^i_j,y^i_j, \bm G_i) - g_i(x^i_j,y^i_j) \right]^2 \right\}, &(x^i_j, y^i_j) \in \bm{X}_{\text{BC}_2,i}. \\
\end{aligned}
\end{equation}

\begin{figure}[ht]
    \centering
    \includegraphics[width=12cm]{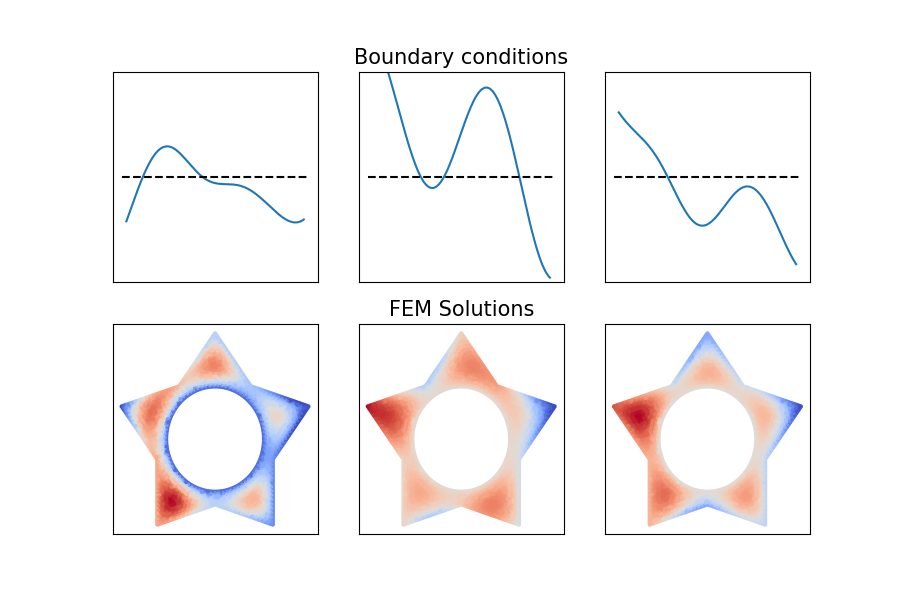}
    \caption{\footnotesize Three representative FEM responses (bottom row) calculated for three realizations of pressure values of the boundary  (top row). The dashed line represents the zero pressure.}
    \label{fig.darcy_samples}
\end{figure}

\textcolor{black}{To better clarify the relationship between the boundary conditions and PDE solutions, some samples of Gaussian process for boundary and PDE solution are shown in Figure \ref{fig.darcy_samples}.} 

\subsubsection{A 2D plate problem}
In this example, we consider more complex problems where a nonlinear mapping between multiple parameter functions and multiple response functions is to be learned. In particular, let us consider a two-dimensional rectangular elastic plate with a central hole (similar to the example in \cite{PINN_plate}). The  steady state solution for the plate displacements is governed by the following system of partial differential equations
\begin{equation}
\begin{aligned}
    \frac{E}{1-\mu^2}\left[\frac{\partial^2 u(x,y)}{\partial x^2} + \frac{(1-\mu)}{2} \frac{\partial^2 u(x,y)}{\partial y^2} + \frac{(1+\mu)}{2} \frac{\partial^2 v(x,y)}{\partial x \partial y}\right] &= 0, \quad (x,y) \in \Omega,  \\
    \frac{E}{1-\mu^2}\left[\frac{\partial^2 v(x,y)}{\partial y^2} + \frac{(1-\mu)}{2} \frac{\partial^2 v(x,y)}{\partial x^2} + \frac{(1+\mu)}{2} \frac{\partial^2 u(x,y)}{\partial x \partial y}\right] &= 0, \quad (x,y) \in \Omega,
\end{aligned}
\end{equation}
where $u$ and $v$ are the plate displacements in $x$ and $y$ directions, respectively, $E$ is the Young's \textcolor{black}{modulus}, and $\mu$ is the Poisson's Ratio. 

\begin{figure}[ht]
    \centering
    \includegraphics[width=0.95\textwidth]{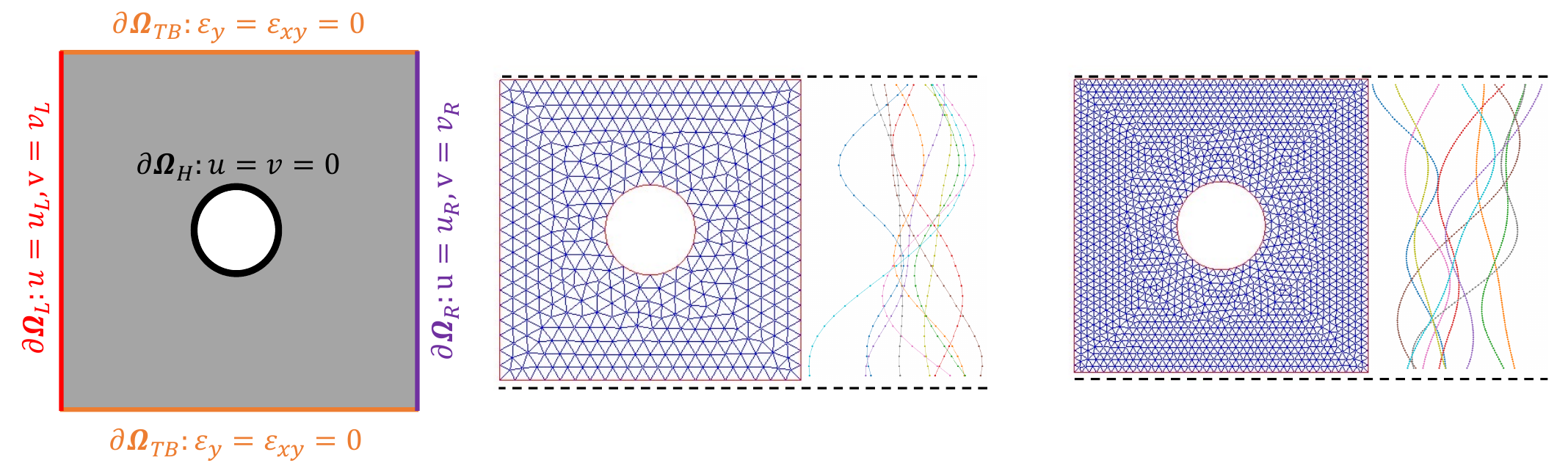}
    \caption{\footnotesize The boundary conditions of 2D plate experiment, and the coarsest and finest meshes used as discretizations to solve the PDE are shown. The samples of Gaussian processes as the prescribed displacement on leftmost edges and rightmost edges are also shown. The dots on the curve are considered as the representation of the prescribed displacements.}
    \label{fig.plate_BC}
\end{figure}

 The  plate is 20 mm × 20 mm, and the hole has a diameter of 5 mm. We impose different boundary conditions on the different edges of the plate. The edge surrounding the central hole is subjected to a fixed boundary condition and there are imposed displacements on the left and right sides. The top and bottom sides are assigned free boundary conditions. The complete boundary conditions are shown in Figure ~\ref{fig.plate_BC}, and are given by
\begin{equation}
\begin{aligned}
    &u(x,y) = 0, \quad v(x,y) = 0, &(x,y) \in \partial \Omega_{\text{H}}, \\
    &u(x,y) = u_l(x,y), \quad v(x,y) = v_l(x,y) & (x,y) \in \partial \Omega_{\text{L}}, \\
    &u(x,y) = u_R(x,y), \quad v(x,y) = v_R(x,y) & (x,y) \in \partial \Omega_{\text{R}}, \\
    &\frac{\partial v(x,y)}{\partial y} = 0, \quad \frac{1}{2} \left[ \frac{\partial u(x,y)}{\partial y} + \frac{\partial v(x,y)}{\partial x} \right] = 0, & (x,y) \in \partial \Omega_{\text{TB}}.
\end{aligned}
\end{equation}
where the  prescribed displacement functions $u_l, v_l, u_R, v_R$ are the inputs to the operator. Similar to previous example, these functions are considered to follow a Gaussian Process model. The covariance kernel  is assumed to be dependent only on the vertical distance, i.e.,
\begin{equation}
\begin{aligned}
    u_l, v_l, u_R, v_R &\sim \mathcal{GP}(1, K(y,y')),\\
    K(y,y') &= \exp\left[-\frac{(y-y')^2}{2l^2}\right]. 
\end{aligned}
\end{equation}
We used $l=5$, and similarly to previous example, we draw $N=1000$ realizations of these functions, and obtain FE solutions using meshes with different sizes.   Figure \ref{fig.plate_BC} shows the finest and coarsest meshes with a few samples of imposed displacement functions. The boundary condition information on the left and right sides  for the $i$-th realization is collected into the following two vectors,
model as follows: 
\begin{equation}
\begin{aligned}
    \bm{X}_{\text{BC}_{\text{H}},i} = \{(x_j^{1,i}, y_j^{1,i})\}_{j=1}^{M_{1,i}}, \\
    \bm{X}_{\text{BC}_{\text{L}},i} = \{(x_j^{2,i}, y_j^{2,i})\}_{j=1}^{M_{2,i}}, \\
    \bm{X}_{\text{BC}_{\text{R}},i} = \{(x_j^{3,i}, y_j^{3,i})\}_{j=1}^{M_{3,i}}, \\
    \bm{X}_{\text{BC}_{\text{TB}},i} = \{(x_j^{4,i}, y_j^{4,i})\}_{j=1}^{M_{4,i}},
\end{aligned}
\end{equation}
where $\bm{X}_{\text{BC}_{\text{H}},i}$, $\bm{X}_{\text{BC}_{\text{L}},i}$, $\bm{X}_{\text{BC}_{\text{R}},i}$, $\bm{X}_{\text{BC}_{\text{TB}},i}$ are the sets of boundary coordinates of $\partial \Omega_{\text{H}}$, $\partial \Omega_{\text{L}}$, $\partial \Omega_{\text{R}}$, $\partial \Omega_{\text{TB}}$, respectively. The prescribed displacement information as inputs to the model $\bm{G}_{L,i}$ and $\bm{G}_{R,i}$ are formulated as:
\begin{equation}
\begin{aligned}
    \bm{G}_{L,i} &= \{\bm{X}_{\text{BC}_L,i}, u_{L,i}(\bm{X}_{\text{BC}_L,i}), v_{L,i}(\bm{X}_{\text{BC}_L,i})\}, \\
    \bm{G}_{R,i} &= \{\bm{X}_{\text{BC}_R,i}, u_{R,i}(\bm{X}_{\text{BC}_R,i}), v_{R,i}(\bm{X}_{\text{BC}_R,i})\}. 
\end{aligned}
\end{equation}

Given this setup, we seek to learn the nonlinear mapping $\mathcal{M}_{\text{Plate}}$ that transforms a given prescribed displacement functions $u_L(x,y), v_L(x,y), u_R(x,y), v_R(x,y)$ on $\partial \Omega_{\text{L}}$ and  $\partial \Omega_{\text{R}}$ to the displacement field $u(x,y)$ and $v(x,y)$ in the entire domain, i.e.,
\begin{equation}
    \mathcal{M}_{\text{Plate}}: [u_L(x,y), v_L(x,y), u_R(x,y), v_R(x,y)]  \rightarrow [u(x,y), v(x,y)].
\end{equation}

Based on our boundary condition settings, the training loss function $\mathcal{L}(\theta)$ for training the neural operators $u_{\theta}(x,y,\bm{G}_L, \bm{G}_R)$ and $v_{\theta}(x,y,\bm{G}_L, \bm{G}_R)$ is formulated as:
\begin{equation}
    \mathcal{L}(\theta) = \alpha_0 \mathcal{L}_{\text{PDE}}(\theta) + \alpha_2 \mathcal{L}_{\text{BC}_{\text{H}}}(\theta) + \alpha_2 \mathcal{L}_{\text{BC}_{\text{L}}}(\theta) + \alpha_3 \mathcal{L}_{\text{BC}_{\text{R}}}(\theta) + \alpha_4 \mathcal{L}_{\text{BC}_{\text{TB}}}(\theta),
\end{equation}
where $\alpha_0$, $\alpha_1$, $\alpha_2$, $\alpha_3$, $\alpha_4$ are the trade-off coefficient, set to be $\alpha_0=0.00001, \alpha_1=\alpha_2=\alpha_3=\alpha_4=1$ to ensure similar orders of magnitude, $\mathcal{L}_{\text{PDE}}(\theta)$ is the PDE residual given by 
\begin{equation}
    \begin{aligned}
    \mathcal{L}_{\text{PDE}}(\theta) = \frac{1}{N} \sum_{i=1}^N \ & \left\{ \frac{E}{1-\mu^2} \left(\frac{\partial^2 u_{\theta}(x,y,\bm{G}_{L,i}, \bm{G}_{R,i})}{\partial x^2} + \frac{\partial^2 v_{\theta}(x,y,\bm{G}_{L,i}, \bm{G}_{R,i})}{\partial y^2} \right) + \right.\\
    & \left. \frac{E}{2(1+\mu)} \left(\frac{\partial^2 u_{\theta}(x,y,\bm{G}_{L,i}, \bm{G}_{R,i})}{\partial y^2} + \frac{\partial^2 v_{\theta}(x,y,\bm{G}_{L,i}, \bm{G}_{R,i})}{\partial x^2} \right) + \right.\\
    & \left.  \frac{E}{2(1-\mu)} \left(\frac{\partial^2 v_{\theta}(x,y,\bm{G}_{L,i}, \bm{G}_{R,i})}{\partial x \partial y} + \frac{\partial^2 u_{\theta}(x,y,\bm{G}_{L,i}, \bm{G}_{R,i})}{\partial x \partial y} \right) \right\}, \quad (x,y) \in \Omega, \\
\end{aligned}
\end{equation}
and $\mathcal{L}_{\text{BC}_{\text{H}}}(\theta)$, $\mathcal{L}_{\text{BC}_{\text{L}}}(\theta)$, $\mathcal{L}_{\text{BC}_{\text{R}}}(\theta)$, $\mathcal{L}_{\text{BC}_{\text{TB}}}(\theta)$ are the BC losses, given by
\begin{equation}
    \begin{aligned}
    \mathcal{L}_{\text{BC}_{\text{H}}}(\theta) =  \frac{1}{N} \sum_{i=1}^N \ & \left\{ \frac{1}{M_{1,i}} \sum_{j=1}^{M_{1,i}}   \left[u_{\theta}(x_j^{i},y_j^{i},\bm{G}_{L,i}, \bm{G}_{R,i}) - 0 \right]^2 + \left[v_{\theta}(x_j^{i},y_j^{i},\bm{G}_{L,i}, \bm{G}_{R,i})- 0 \right]^2 \right\},   \\
    & \qquad \qquad \qquad \qquad \qquad \qquad \qquad \qquad \qquad \qquad \qquad \qquad \qquad \qquad \qquad \qquad {(x_j^{i},y_j^{i}) \in \bm{X}_{\text{BC}_{\text{H}},i},}\\
    \mathcal{L}_{\text{BC}_{\text{L}}}(\theta) = \frac{1}{N} \sum_{i=1}^N \ & \left\{ \frac{1}{M_{2,i}} \sum_{j=1}^{M_{2,i}}   \left[u_{\theta}(x_j^{i},y_j^{i},\bm{G}_{L,i}, \bm{G}_{R,i}) - u_{L,i}(x_j^{i},y_j^{i}) \right]^2 + \left[v_{\theta}(x_j^{i},y_j^{i},\bm{G}_{L,i}, \bm{G}_{R,i})- v_{L,i}(x_j^{i},y_j^{i}) \right]^2 \right\},  \\
    & \qquad \qquad \qquad \qquad \qquad \qquad \qquad \qquad \qquad \qquad \qquad \qquad \qquad \qquad \qquad \qquad {(x_j^{i},y_j^{i}) \in \bm{X}_{\text{BC}_{\text{L}},i},}\\
    \mathcal{L}_{\text{BC}_{\text{R}}}(\theta) = \frac{1}{N} \sum_{i=1}^N \ & \left\{ \frac{1}{M_{3,i}} \sum_{j=1}^{M_{3,i}}   \left[u_{\theta}(x_j^{i},y_j^{i},\bm{G}_{L,i}, \bm{G}_{R,i}) - u_{R,i}(x_j^{i},y_j^{i}) \right]^2 + \left[v_{\theta}(x_j^{i},y_j^{i},\bm{G}_{L,i}, \bm{G}_{R,i})- v_{R,i}(x_j^{i},y_j^{i}) \right]^2 \right\},  \\
    & \qquad \qquad \qquad \qquad \qquad \qquad \qquad \qquad \qquad \qquad \qquad \qquad \qquad \qquad \qquad \qquad {(x_j^{i},y_j^{i}) \in \bm{X}_{\text{BC}_{\text{R}},i},}\\
    \mathcal{L}_{\text{BC}_{\text{TB}}}(\theta) =  \frac{1}{N} \sum_{i=1}^N \ & \left\{ \frac{1}{M_{4,i}} \sum_{j=1}^{M_{4,i}}   \left[ \frac{\partial v_{\theta}(x_j^{i},y_j^{i},\bm{G}_{L,i}, \bm{G}_{R,i})}{\partial y_j^{i}} - 0 \right]^2 + \right.  \\ 
    & \left. \left[ \frac{1}{2} \left( \frac{\partial u_{\theta}(x_j^{i},y_j^{i},\bm{G}_{L,i}, \bm{G}_{R,i})}{\partial y_j^{i}} + \frac{\partial v_{\theta}(x_j^{i},y_j^{i},\bm{G}_{L,i}, \bm{G}_{R,i})}{\partial x_j^{i}} \right) - 0\right]^2 \right\},  \\
    & \qquad \qquad \qquad \qquad \qquad \qquad \qquad \qquad \qquad \qquad \qquad \qquad \qquad \qquad \qquad \qquad {(x_j^{i},y_j^{i}) \in \bm{X}_{\text{BC}_{\text{TB}},i}.}
\end{aligned}
\end{equation}

\begin{figure}[ht]
    \centering
    \includegraphics[width=16cm]{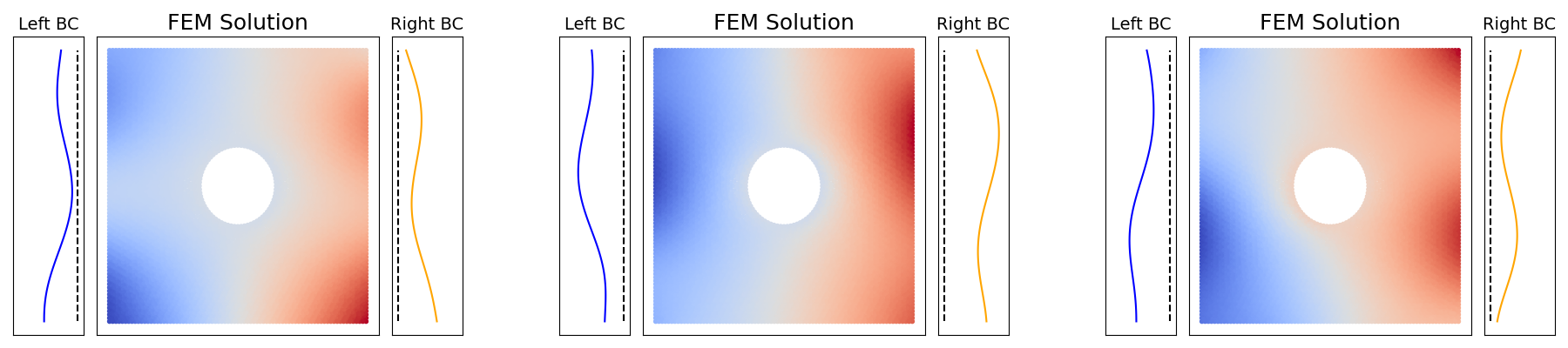}
    \caption{\footnotesize Three representative samples of the FEM responses with their corresponding   prescribed displacements on the left  and right edges. The dashed line represents  zero displacement.}
    \label{fig.plate_samples}
\end{figure}

\textcolor{black}{To better clarify the relationship between the boundary conditions and PDE solutions, some samples of Gaussian process for boundary and PDE solution are shown in Figure \ref{fig.plate_samples}.} 

\subsection{Main results}
For each problem, we conduct two experiments. First, in order to compare our method with PI-DeepONet, which can only handle data with the same \textcolor{black}{ discrete representation},  we consider a fixed \textcolor{black}{discretization size}, taken to be small, for all the training and test cases. In the second experiment, we relax this constraint and generate training and test data each with variable \textcolor{black}{ discretization sizes}. 
The comparison between our proposed methods and the PI-DeepONet is summarized in Table ~\ref{fig.PI_com}.  
We report the average error calculated across the entire test dataset, along with the standard deviation of the relative errors, to demonstrate the model stability. Our proposed model showed satisfactory mean performance in both the Darcy flow problem and the 2D plate problem. \textcolor{black}{We found out the improved DeepONet architecture of \cite{IDON} to be the best baseline model. Compared to this baseline architecture, our model achieved 19.4\% accuracy improvement for the Darcy flow problem and 11.9\% accuracy improvement for the 2D plate problem. When the standard deviation of relative errors are considered, our model also shows higher stability.}

In the second experiment, when variable \textcolor{black}{variable discrete representations are used}, compared to the \textcolor{black}{fixed discretization case}, we note that PI-DCON offers smaller accuracy levels for both the Darcy flow and the 2D plate problems. These results are reasonable, especially considering the 2D plate problem is inherently more complex. Despite these challenges, our model is able to maintain acceptable levels of accuracy, with mean relative errors below 3.5\% for both problems.

\begin{table*}[!ht]
\begin{center}
\caption{\textcolor{black}{Accuracy comparison between PI-DCON and PI-DeepONet models. In the experiments with variable discrete representation, we show the model performance across different meshes.} }
\begin{tabular}{c c c c}
\hline
\multirow{2}{*}{Experiment} & \multicolumn{3}{c}{\makecell{$L_2$ relative error}}\\
 {} & PI-DeepONet & \cite{IDON} & PI-DCON \\
\hline
\makecell{Darcy flow \\ (fixed discrete representation)} & 7.10\% $\pm$ 4.00\% & 3.10\% $\pm$ 1.41\% & \bf{2.50\%} $\pm$ \bf{1.11\%}\\
\makecell{2D plate \\ (fixed discrete representation)} & 5.58\% $\pm$ 1.14\% & 2.01\% $\pm$ 0.53\% & \bf{1.77\%} $\pm$ \bf{0.42\%} \\
\hline
\makecell{Darcy flow \\ (variable discrete representation)} & - & - & 3.42\% $\pm$ 1.43\%\\
\makecell{2D plate \\ (variable discrete representation)} & - & - & 2.98\% $\pm$ 0.92\% \\
\hline
\label{fig.PI_com}
\end{tabular}
\end{center}
\end{table*}

\begin{figure}[ht]
    \centering
            \begin{subfigure}[b]{0.3\textwidth}
            \centering
            \includegraphics[width=\textwidth]{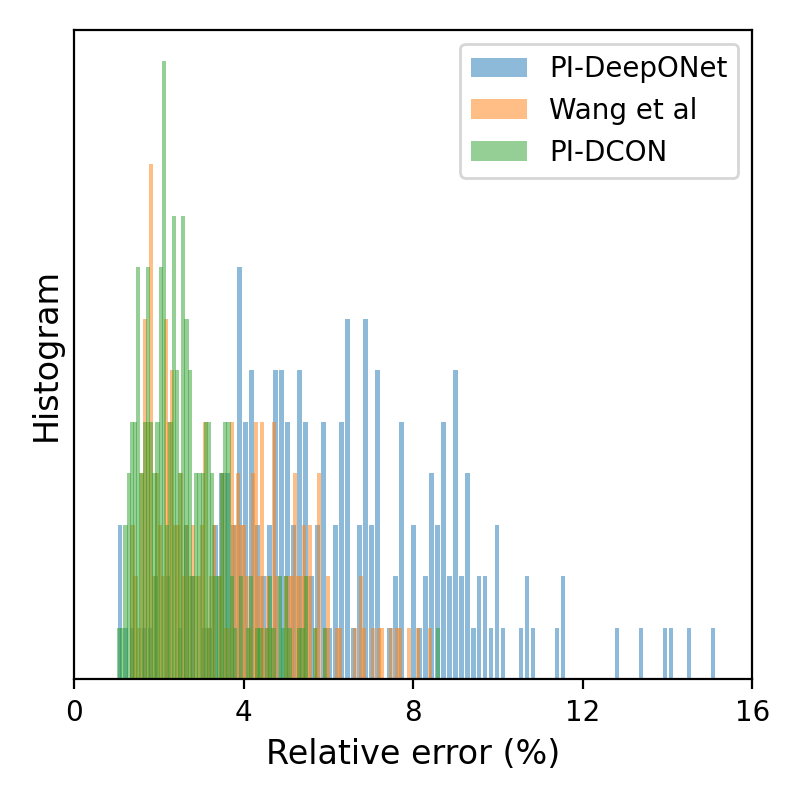}
            \caption{Darcy flow}    
            \label{fig:DCON_vs_DON_darcy_fixed}
        \end{subfigure}
        \hspace{0.2cm}
        \begin{subfigure}[b]{0.3\textwidth}   
            \centering 
            \includegraphics[width=\textwidth]{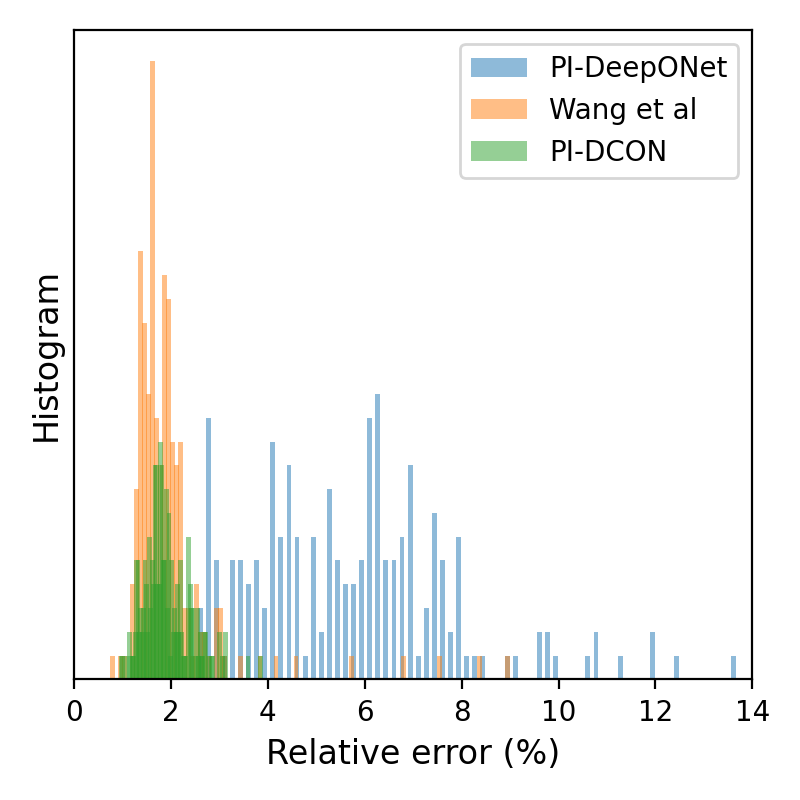}
            \caption{2D plate}  
            \label{fig:DCON_vs_DON_plate_fixed}
        \end{subfigure}
    \caption{\footnotesize Comparison between PI-DCON and PI-DeepONet  on the Darcy flow and 2D plate problems, with the fixed \textcolor{black}{discretization size}.} \label{fig.err_dist.fixed}
\end{figure}

\begin{figure}[ht]
    \centering
        \begin{subfigure}[b]{0.3\textwidth}  
            \centering 
            \includegraphics[width=\textwidth]{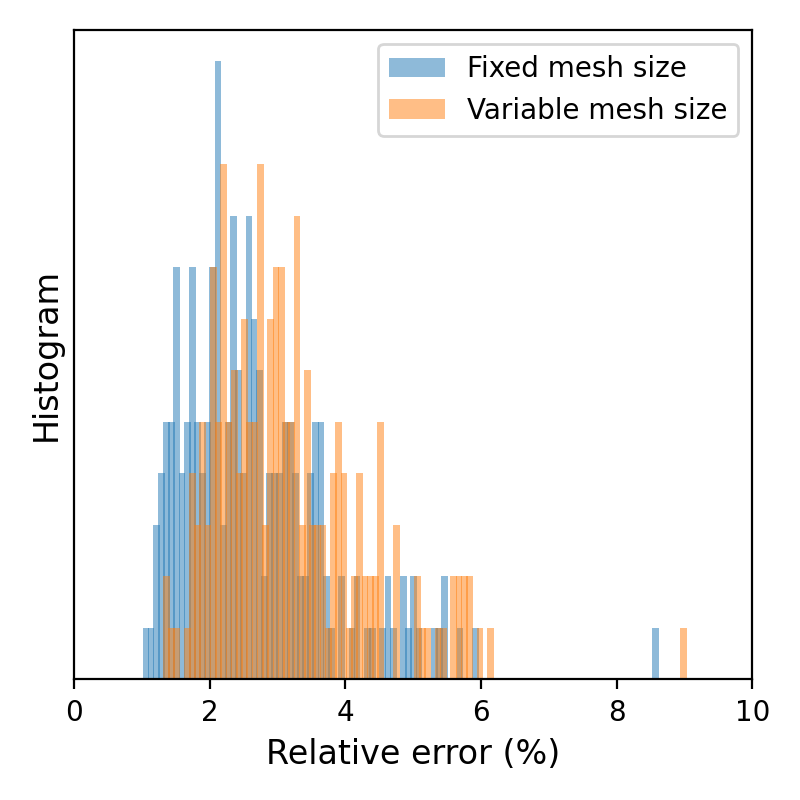}
            \caption{Darcy flow}    
            \label{fig:DCON_darcy_fixed_vs_darcy_vary}
        \end{subfigure}
        \hspace{0.2cm}
        \begin{subfigure}[b]{0.3\textwidth}   
            \centering 
            \includegraphics[width=\textwidth]{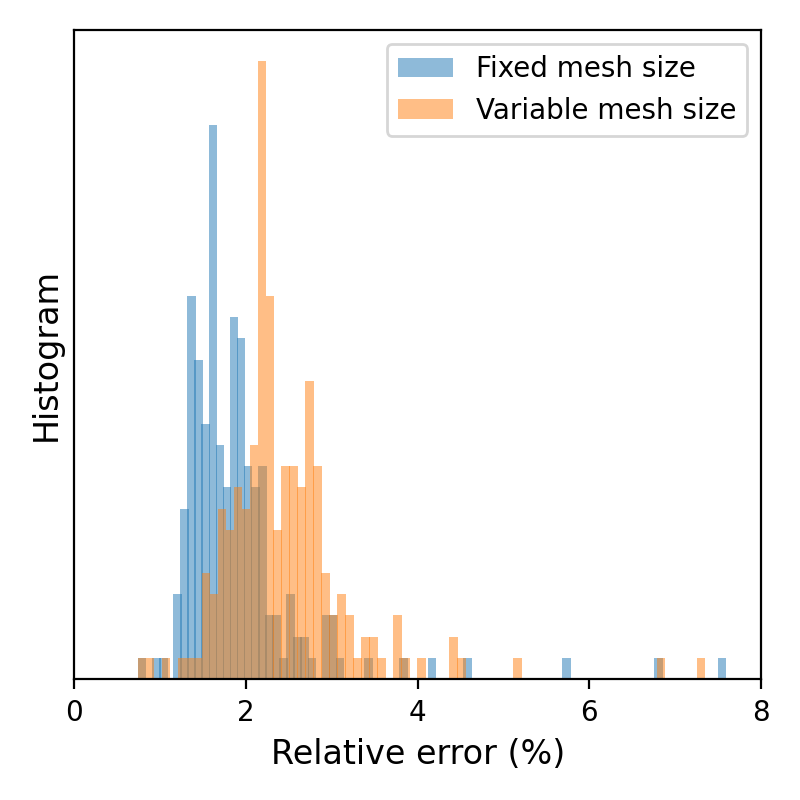}
            \caption{2D plate}  
            \label{fig:DCON_plate_fixed_vs_plate_vary}
        \end{subfigure}
    \caption{\footnotesize Comparison of PI-DCON for fixed \textcolor{black}{discretization size versus variable discretization size}, on the Darcy flow and 2D plate problems. \label{fig.err_dist.var}}
\end{figure}

In order to better visualize our model performances, we present histograms of the estimated relative error across the entire test dataset in Figures~\ref{fig.err_dist.fixed} and \ref{fig.err_dist.var}. Figure~\ref{fig.err_dist.fixed} demonstrates that, for most instances, our model predicts the PDE solutions with a relative error under 4\%, whereas the relative error for most of the predictions of PI-DeepONet falls within the 4\% to 10\% range. \textcolor{black}{Compared with \cite{IDON}, our model also show slightly higher performance.} In Figure~\ref{fig.err_dist.var}, we also show that our model can handle different \textcolor{black}{discretization sizes} with similar prediction accuracy as the fixed \textcolor{black}{discretization size}, achieving relative error under 6\% for predictions of most cases. Figure~\ref{fig.err_dist.var} further illustrates that our model maintains comparable prediction accuracy across different \textcolor{black}{discretization sizes}. For the majority of predictions, it manages to keep the relative error below 6\%.

Additionally, it's noted that each model occasionally produces significantly poor predictions for certain test dataset samples. These instances are considered as outliers within the sample distribution. This phenomenon suggests that all the trained models demonstrate a deficiency in sustaining acceptable performance on samples that fall outside of the distribution, which will be an interesting problem to be investigated in the future works.

We also investigate the average model performance over the whole domain in the geometry in Table ~\ref{fig.avg_performance}. Overall, our model achieved better performance than PI-DeepONet. Also, We observe that the models have difficulties in predicting the PDE solutions on the boundaries of varying boundary conditions. This indicates a larger trade-off coefficient for the hard constraints of the varying boundary conditions may be needed. However, setting this coefficient too high could impede the minimization of PDE residuals, indicating the need for future research to find an optimal balance for the trade-off coefficient.

When examining errors along the boundaries, we particularly notice that larger errors tend to appear at the boundary's sharp angles, such as the concave angles in the pentagram and the corner angles of the 2D plate. This discrepancy could stem from imprecise gradient computations at these points, highlighting the importance of paying extra attention to the sharp angles in the geometry during physics-informed training processes.

\begin{table*}[!ht]
\centering 
\caption{Comparison between  Mean Absolute Errors shown over the domain}
\label{fig.best_worst_cases_darcy} 
\begin{tabular}{|c|c c c|}
\hline
{ } & \makecell{PI-DeepONet \\(fixed discrete representation) }& \makecell{PI-DCON \\ (fixed discrete representation)} & \makecell{PI-DCON \\ (variable discrete representation)} \\ 
\hline
 \rotatebox[origin=l]{90}{\quad Darcy flow} & {\includegraphics[width=0.25\textwidth]{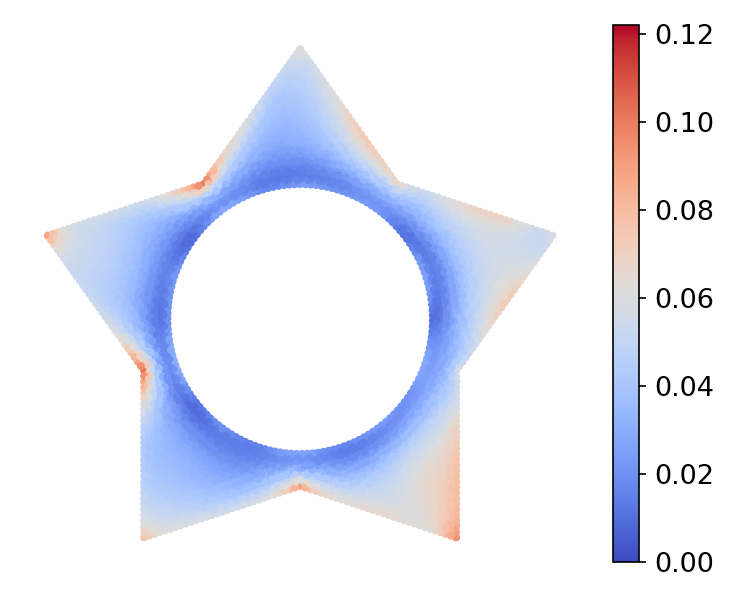}} & {\includegraphics[width=0.25\textwidth]{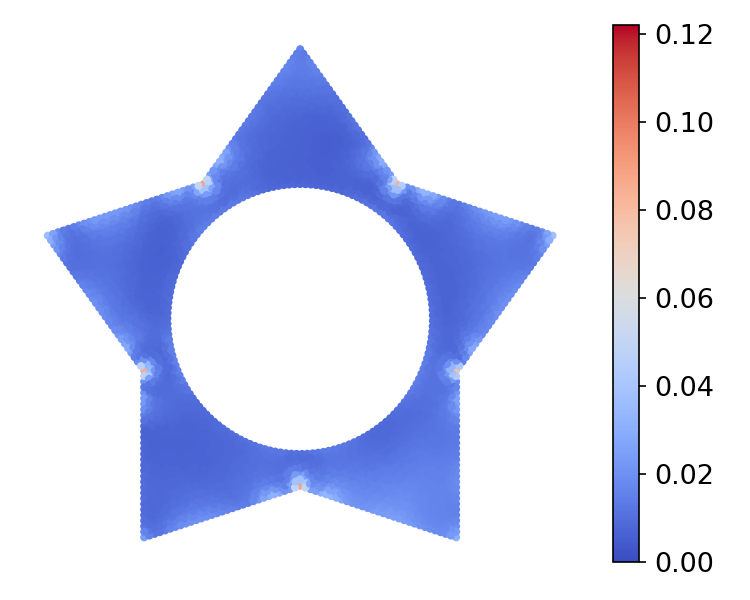}} & {\includegraphics[width=0.25\textwidth]{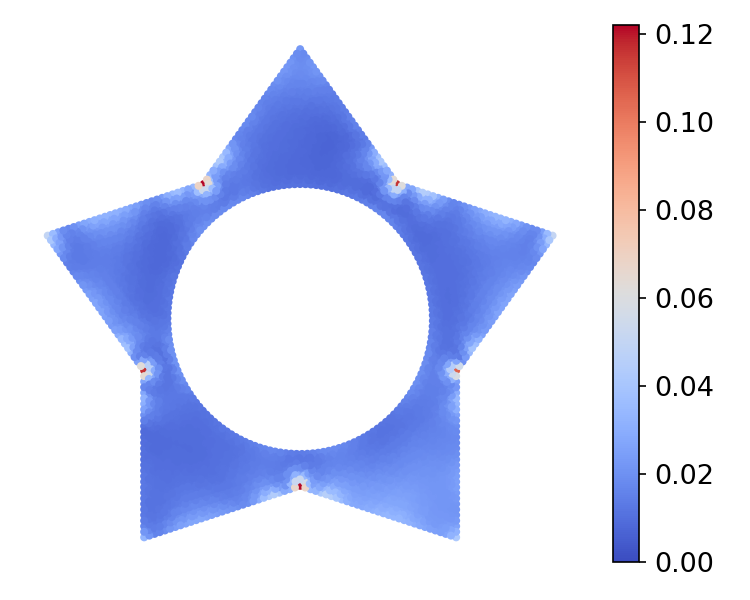}} \\
 \rotatebox[origin=l]{90}{\qquad 2D plate} & {\includegraphics[width=0.25\textwidth]{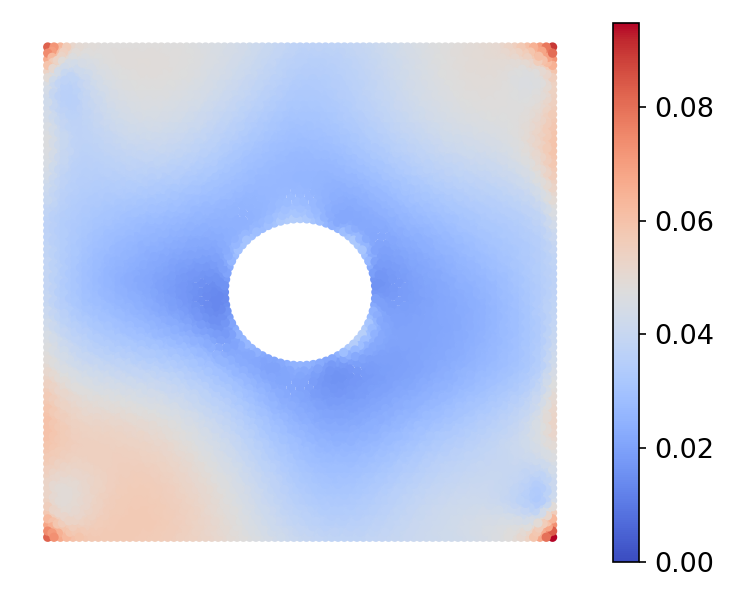}} & {\includegraphics[width=0.25\textwidth]{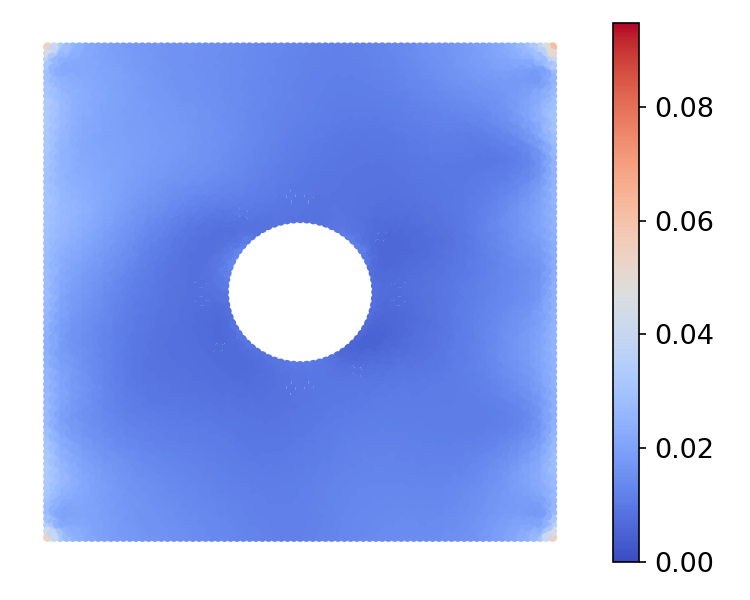}} & {\includegraphics[width=0.25\textwidth]{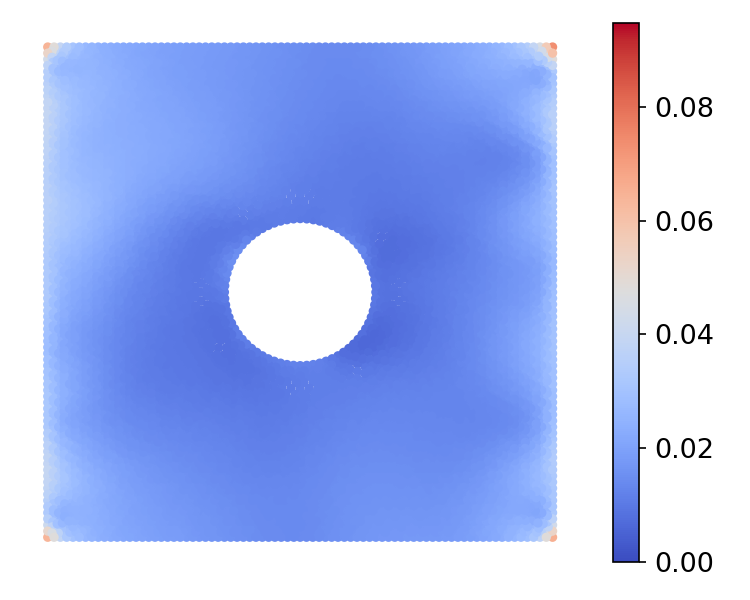}} \\ 
\hline
\end{tabular}
\label{fig.avg_performance}
\end{table*}

To further illustrate the comparative performance of our model, Tables~ \ref{fig.best_worst_case}    show predicted and ground truth solutions for PI-DeepONet and PI-DCON.  For each problem, four of the test cases are shown. In particular, the test cases (i.e. realizations of boundary conditions) that cause  the best and worst performances  of PI-DeepONet and PI-DCON. It can be seen that our model shows higher accuracy, and  maintains remarkable accuracy even in the most challenging (worst-case) scenario. Furthermore, the relatively small performance gap between the worst and best cases also highlights the robustness of our model. Furthermore, as can be seen in Table~\ref{fig.best_worst_case_vary_mesh},  PI-DCON can also produce acceptable predictions when training and test cases are created with variable \textcolor{black}{discretization sizes}.

\begin{table*}[htbp]
\centering 
\caption{Comparison between the performance of PI-DCON and PI-DeepONet. Out of all realizations of boundary conditions, the ones that causes best and worst performances of each model are shown.}
\begin{tabular}{|c|c|c | c c | c c|}
\hline
\multicolumn{2}{|c|}{ } & \multirow{2}{*}{Ground Truth} & \multicolumn{2}{c|}{PI-DeepONet} & \multicolumn{2}{c|}{PI-DCON} \\ 
\cline{4-7}
\multicolumn{2}{|c|}{ } &  & Prediction & Absolute Error & Prediction & Absolute Error \\ 
\hline
\multirow{8}{*}{\rotatebox[origin=c]{90}{Darcy flow }} & \rotatebox[origin=l]{90}{\makecell{PI-DeepONet \\ Best case}} & {\includegraphics[width=0.16\textwidth]{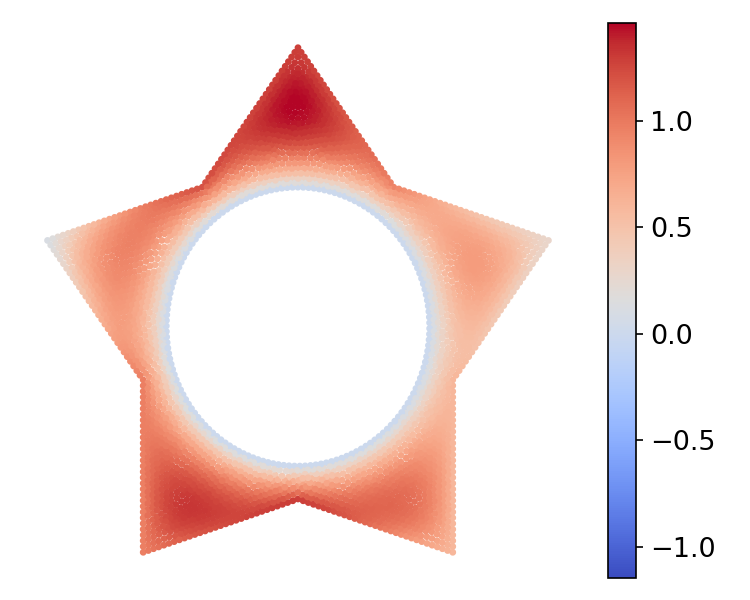}} & {\includegraphics[width=0.16\textwidth]{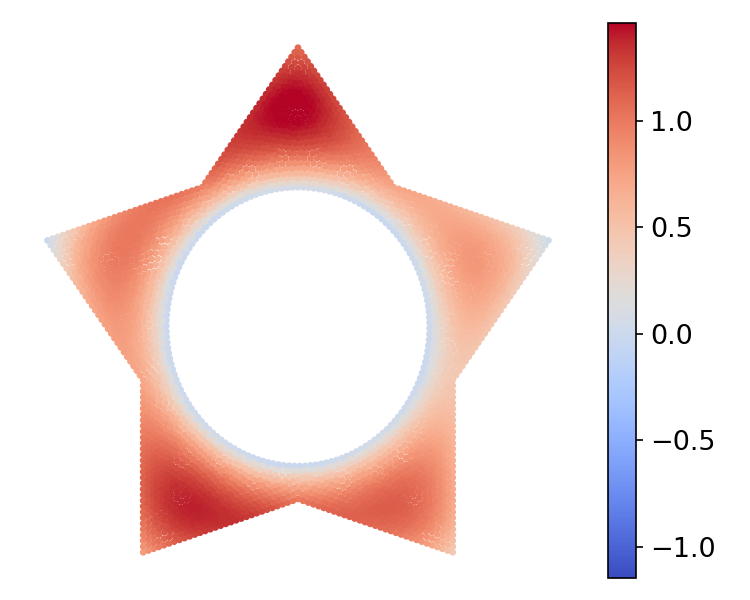}} & {\includegraphics[width=0.16\textwidth]{ 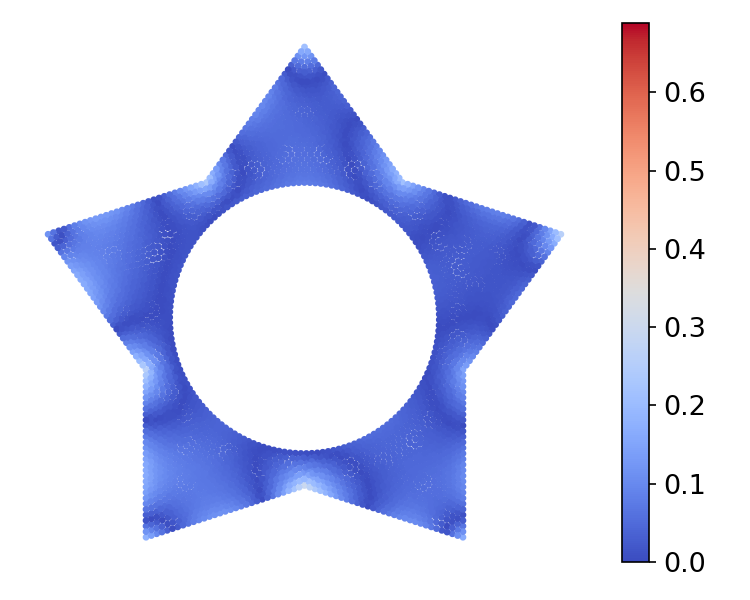}} & {\includegraphics[width=0.16\textwidth]{ 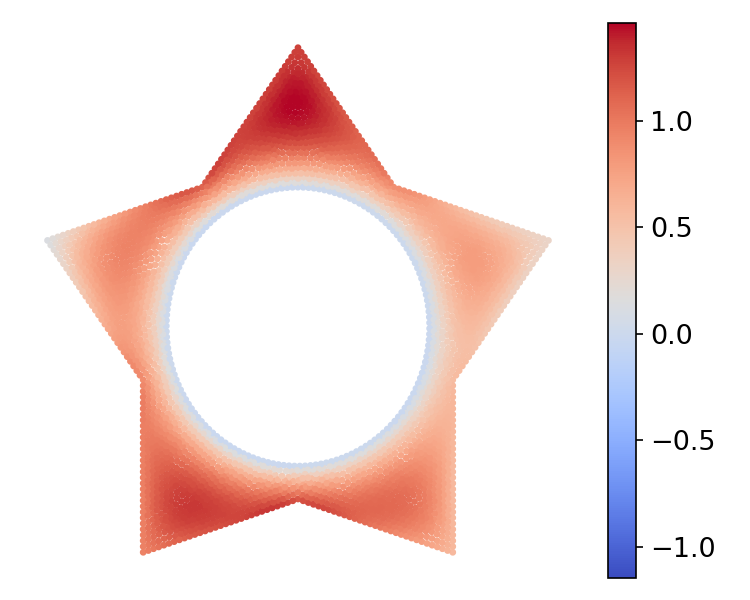}} & {\includegraphics[width=0.16\textwidth]{ 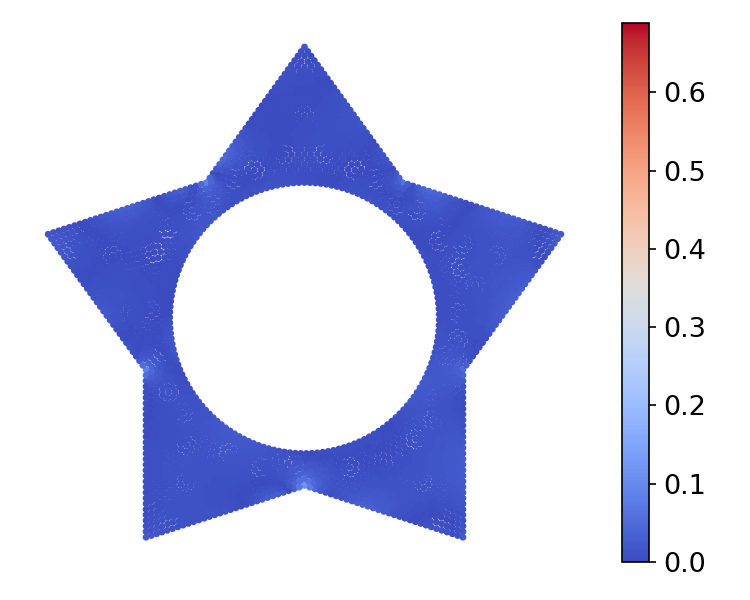}} \\
& \rotatebox[origin=l]{90}{\makecell{PI-DCON \\ Best case}} & {\includegraphics[width=0.16\textwidth]{ 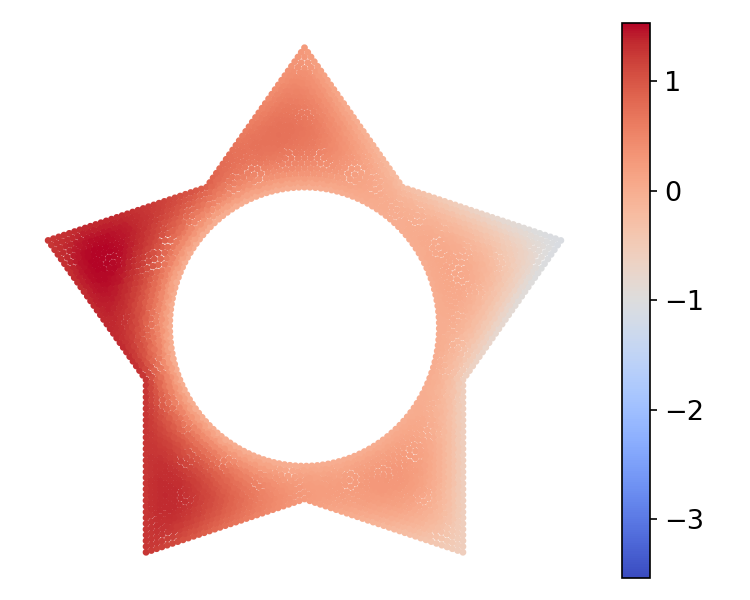}} & 
 {\includegraphics[width=0.16\textwidth]{ 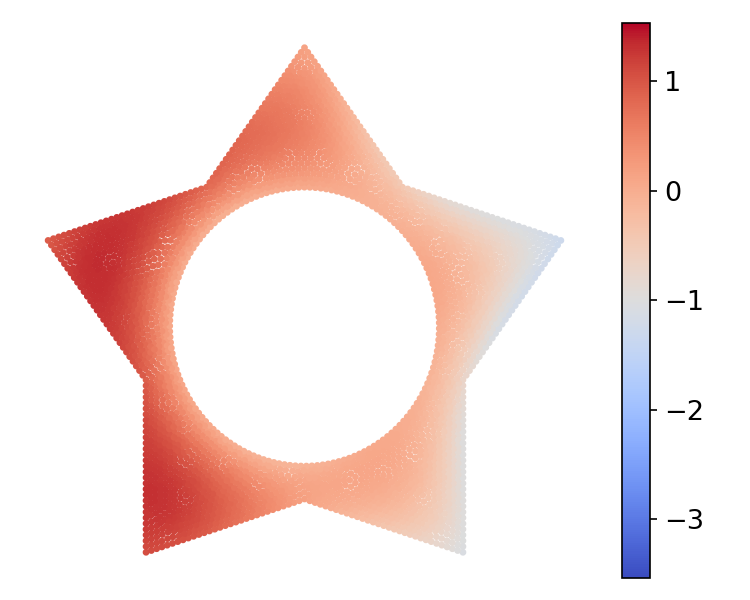}} & {\includegraphics[width=0.16\textwidth]{ 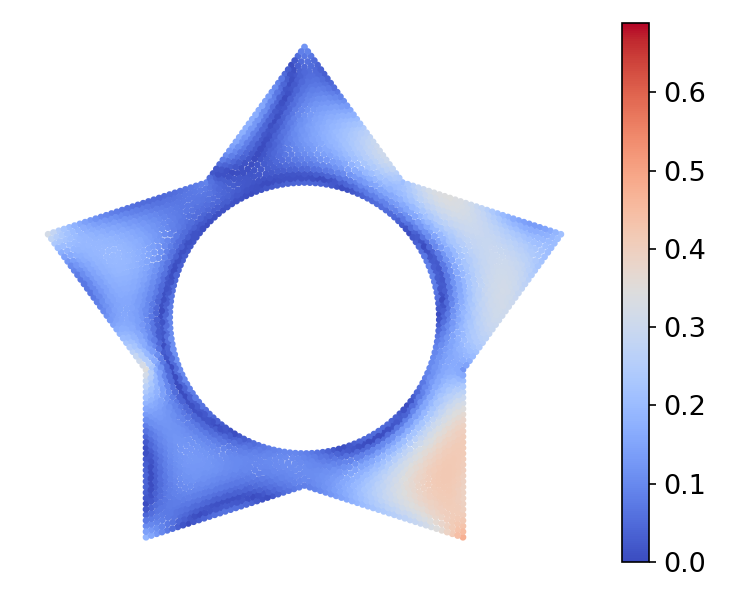}} &{\includegraphics[width=0.16\textwidth] { 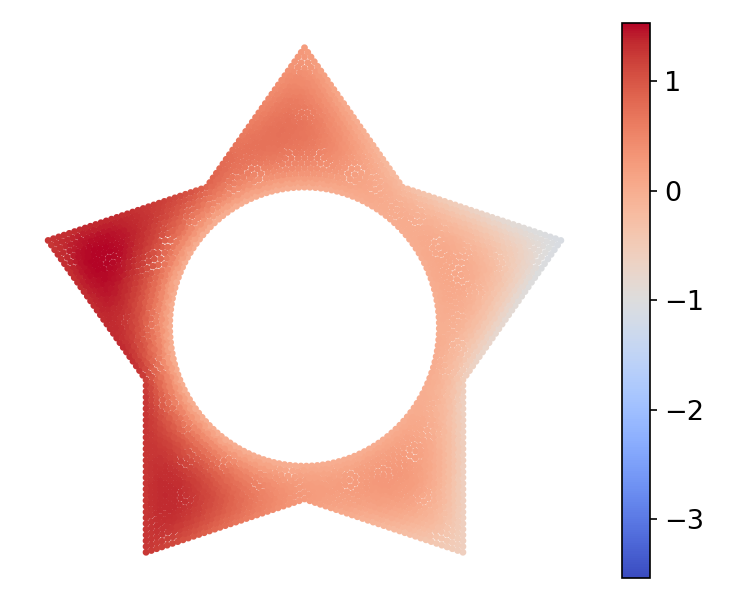}} & {\includegraphics[width=0.16\textwidth]{ 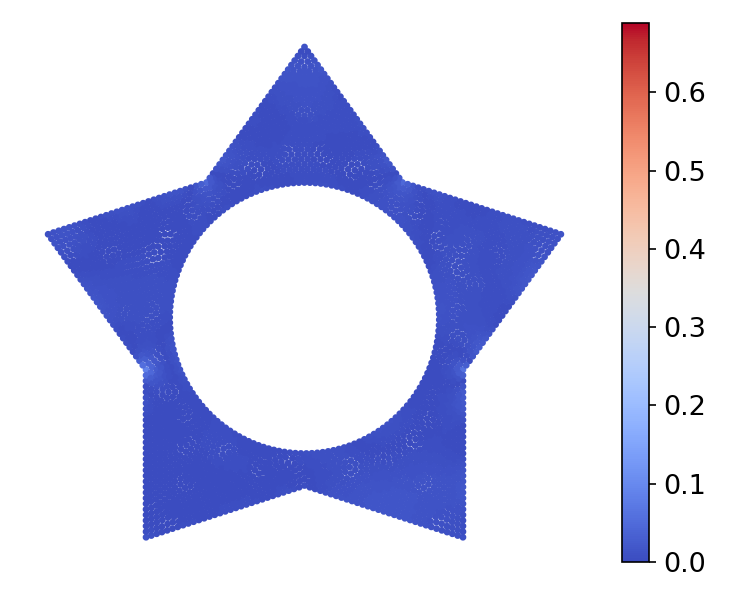}}  \\  
& \rotatebox[origin=l]{90}{\makecell{PI-DeepONet \\ Worst case}} & {\includegraphics[width=0.16\textwidth]{ 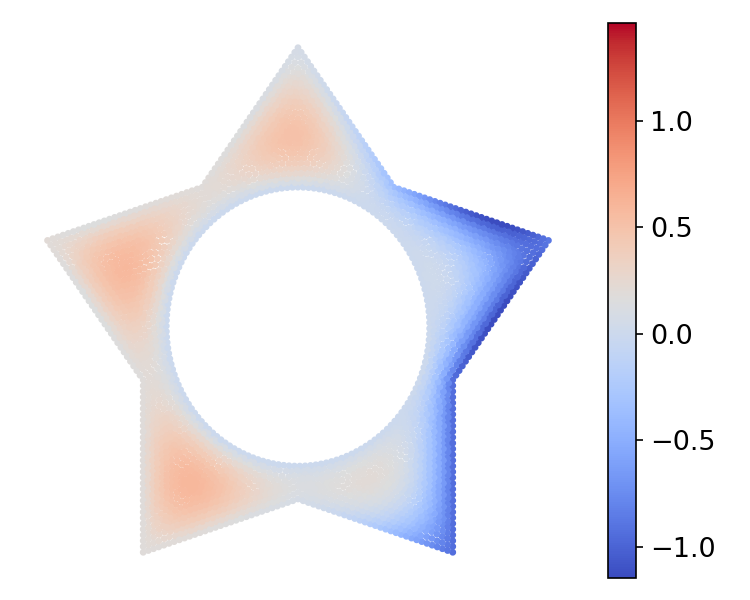}} & {\includegraphics[width=0.16\textwidth]{ 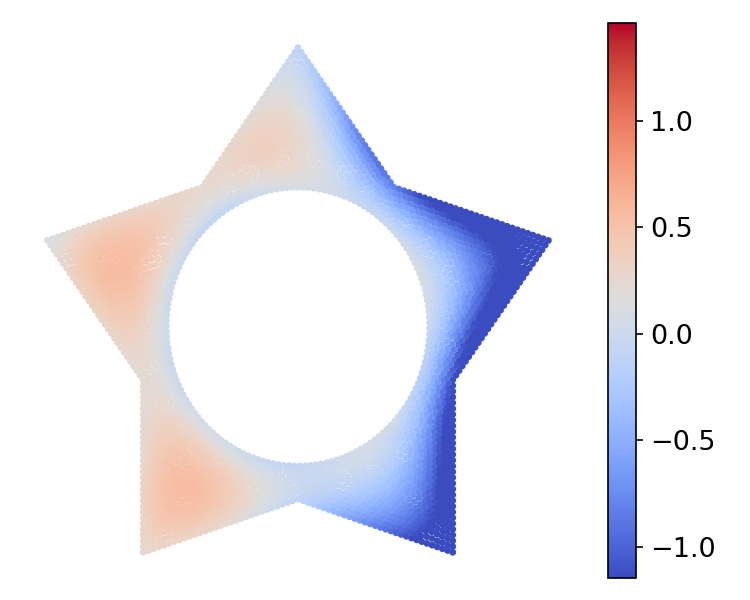}} & {\includegraphics[width=0.16\textwidth]{ 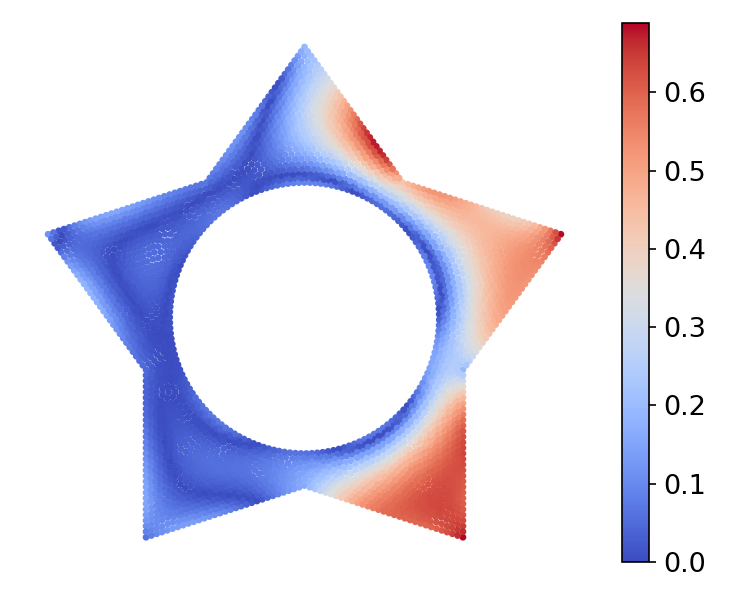}} & {\includegraphics[width=0.16\textwidth]{ 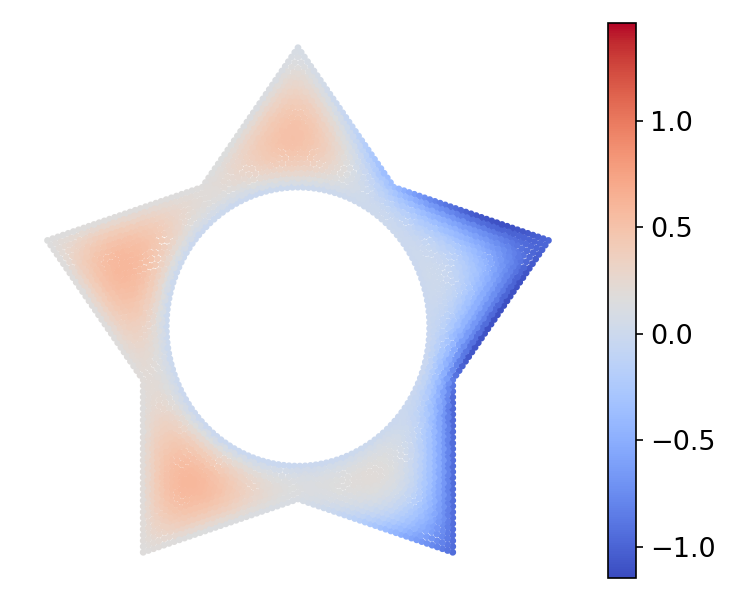}} & {\includegraphics[width=0.16\textwidth]{ 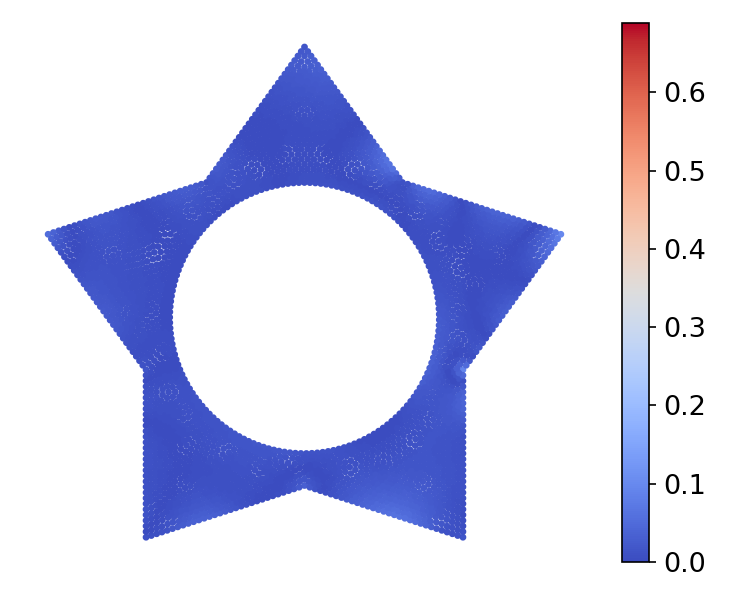}} \\
& \rotatebox[origin=l]{90}{\makecell{PI-DCON \\ Worst case}} & {\includegraphics[width=0.16\textwidth]{ 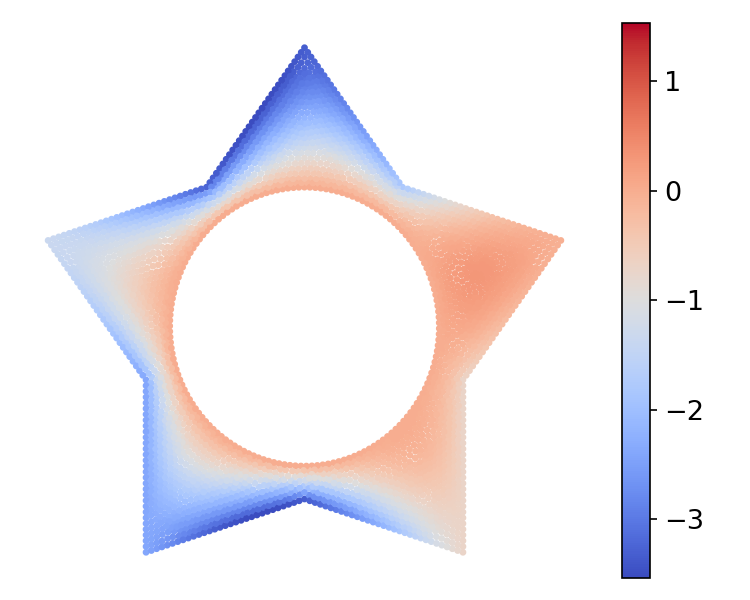}} &
 {\includegraphics[width=0.16\textwidth]{ 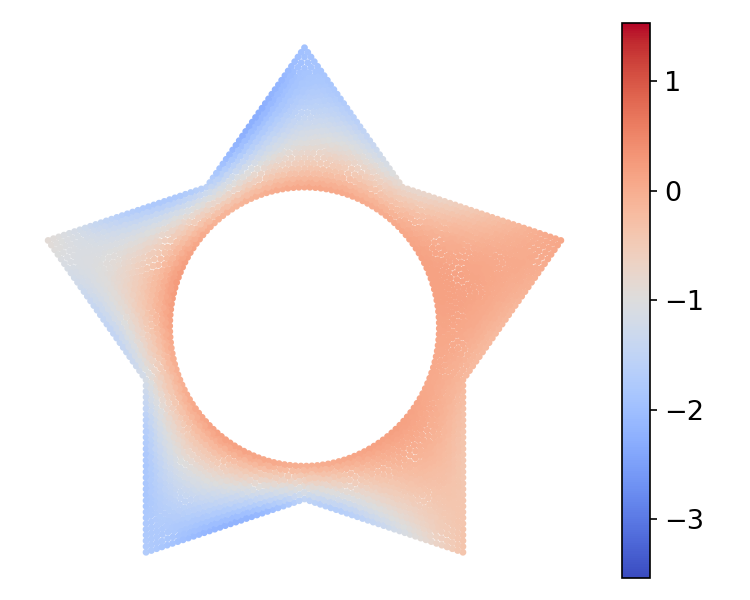}} & {\includegraphics[width=0.16\textwidth]{ 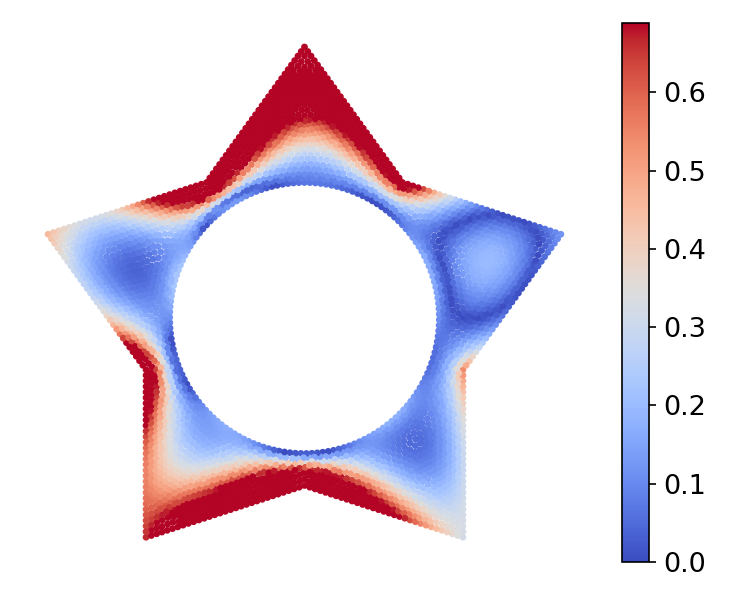}} &{\includegraphics[width=0.16\textwidth]{ 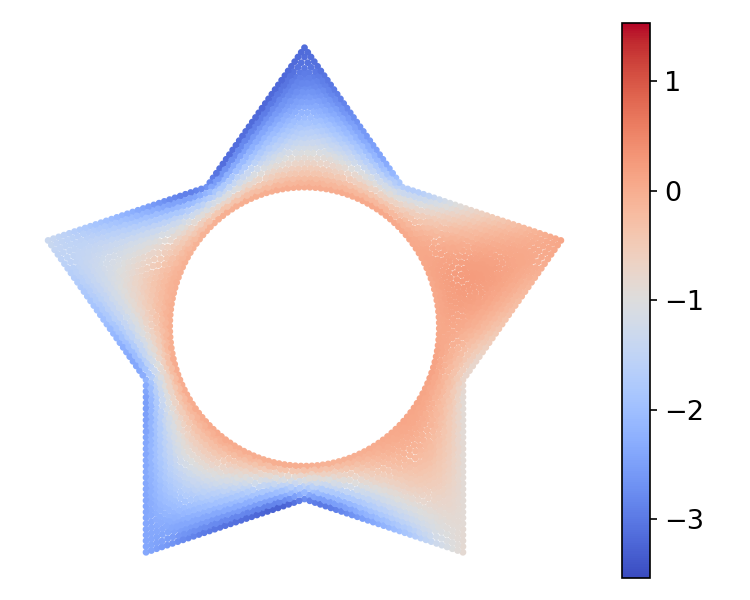}} & {\includegraphics[width=0.16\textwidth]{ 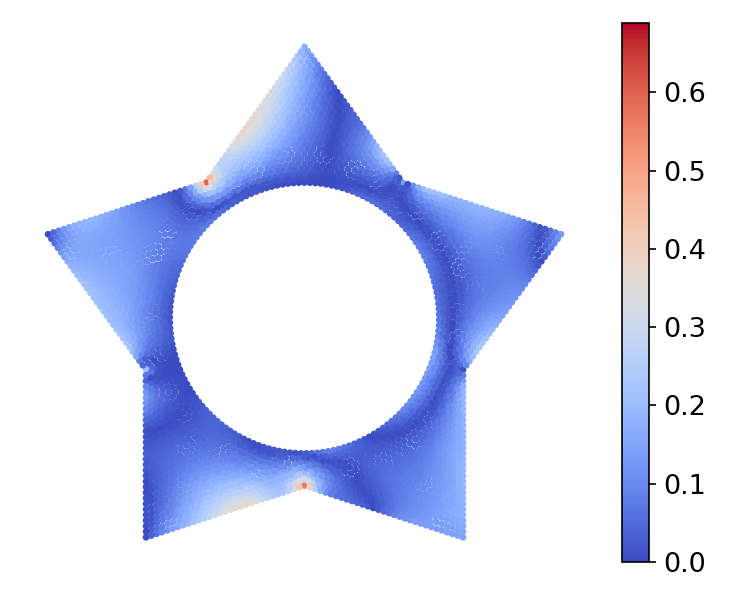}}  \\ 
\hline
\multirow{8}{*}{\rotatebox[origin=c]{90}{2D Plate}} & \rotatebox[origin=l]{90}{\makecell{PI-DeepONet \\ Best case}} & {\includegraphics[width=0.16\textwidth]{ 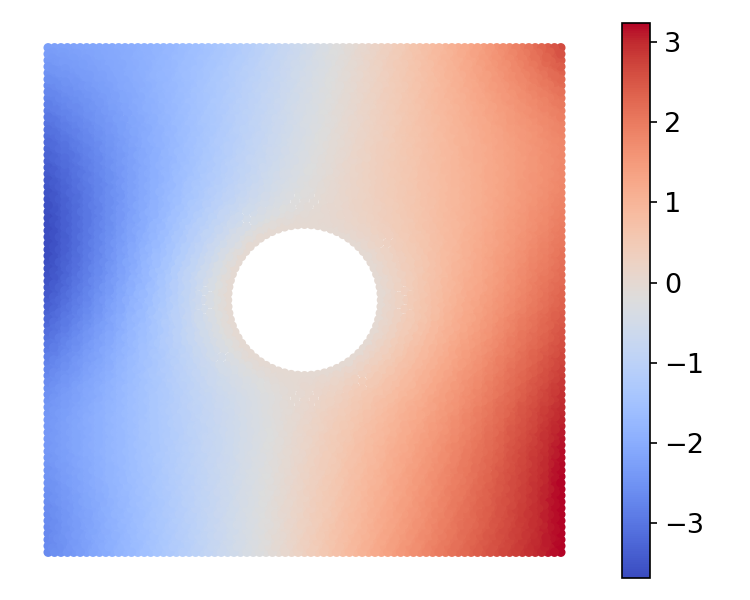}} & {\includegraphics[width=0.16\textwidth]{ 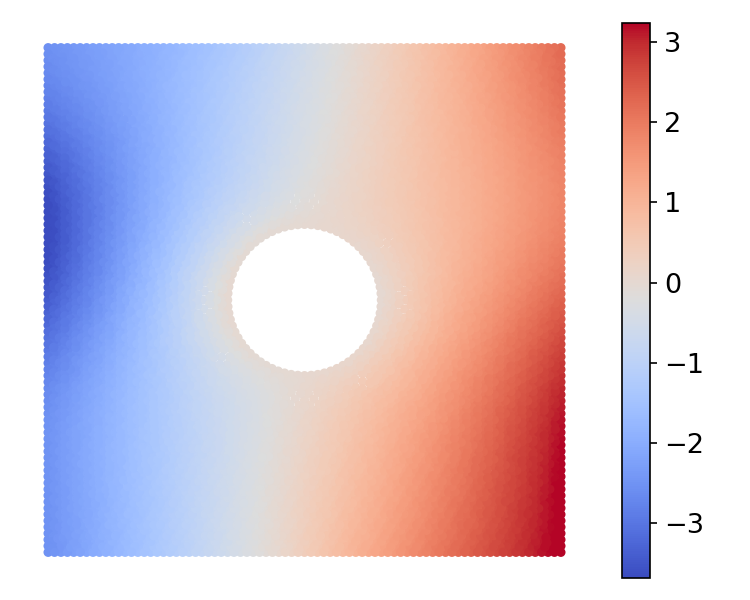}} & {\includegraphics[width=0.16\textwidth]{ 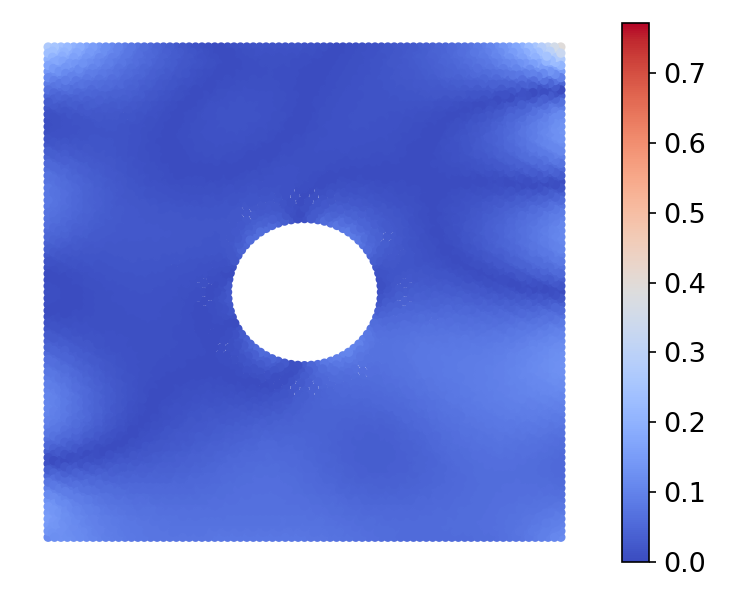}} & 
{\includegraphics[width=0.16\textwidth]{ 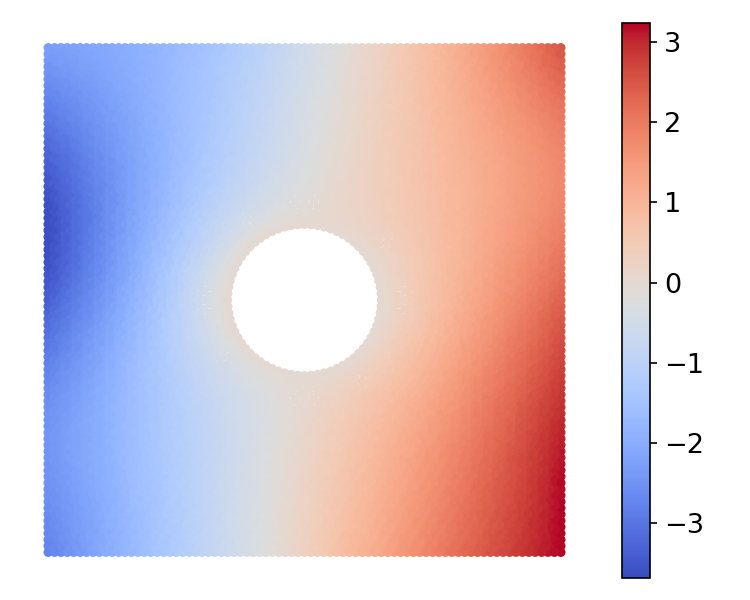}} & {\includegraphics[width=0.16\textwidth]{ 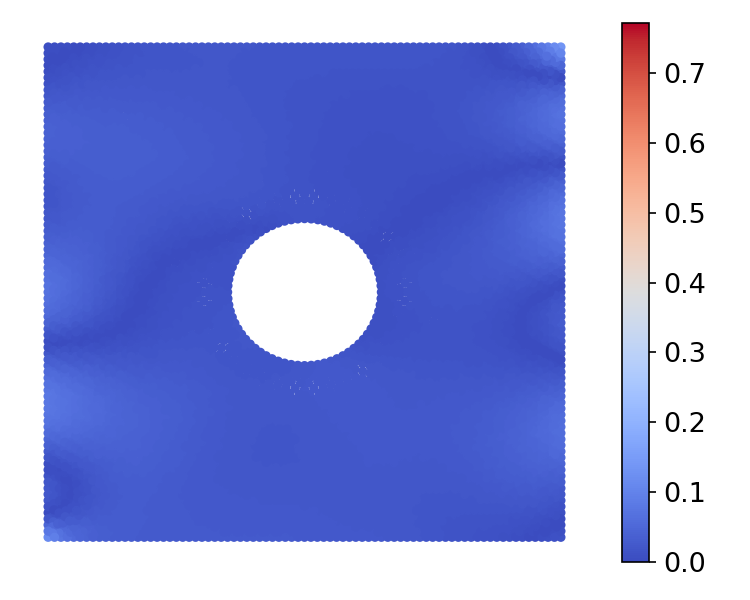}}\\
 & \rotatebox[origin=l]{90}{\makecell{PI-DCON \\ Best case}} & {\includegraphics[width=0.16\textwidth]{ 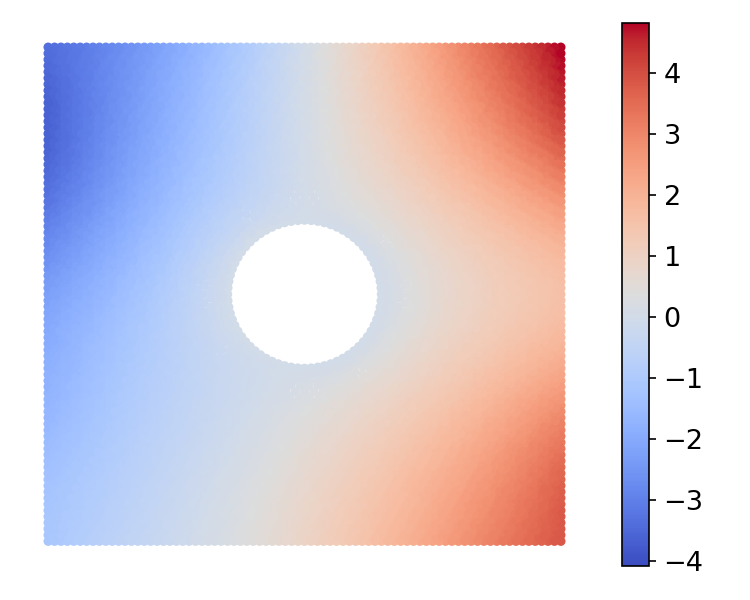}} & 
 {\includegraphics[width=0.16\textwidth]{ 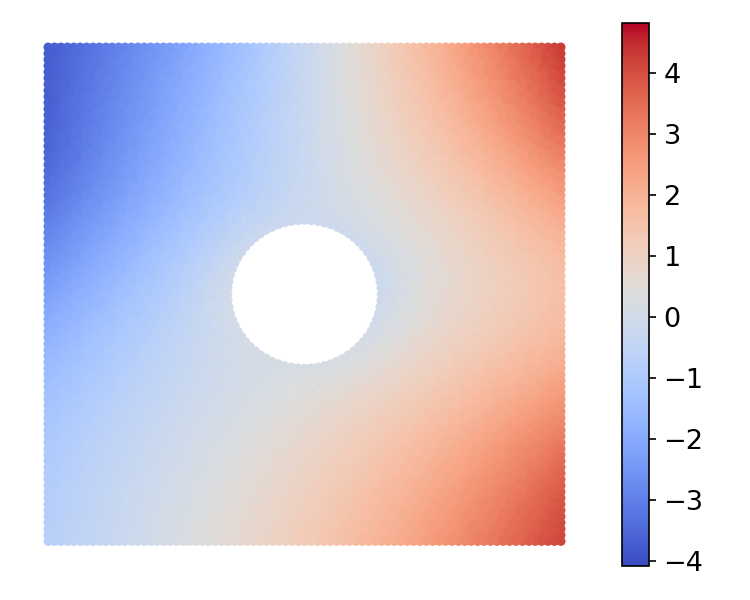}} & {\includegraphics[width=0.16\textwidth]{ 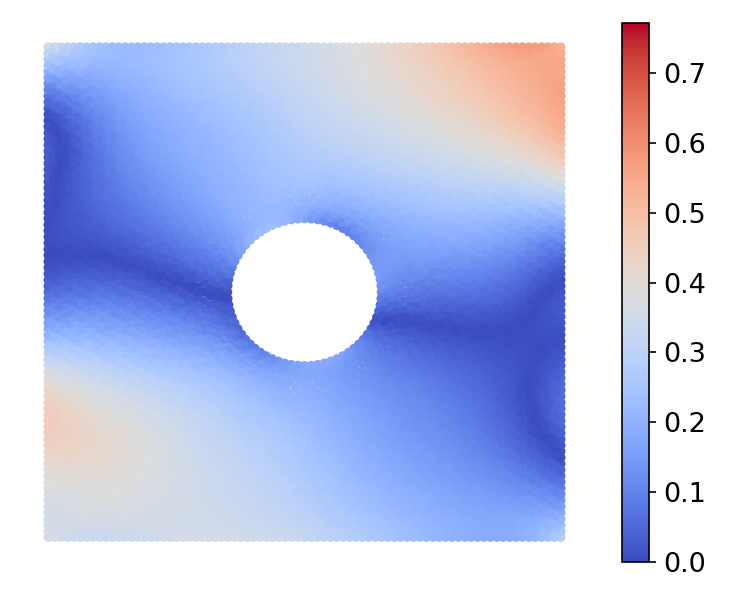}} &{\includegraphics[width=0.16\textwidth]{ 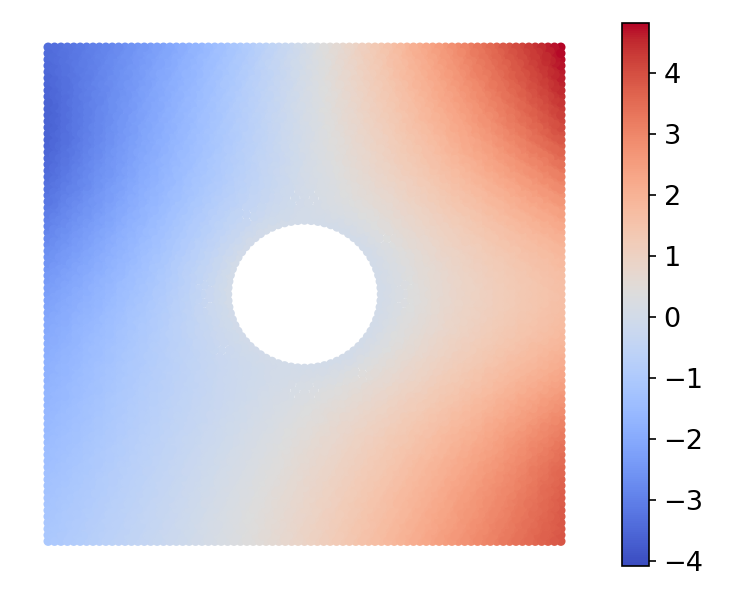}} & {\includegraphics[width=0.16\textwidth]{ 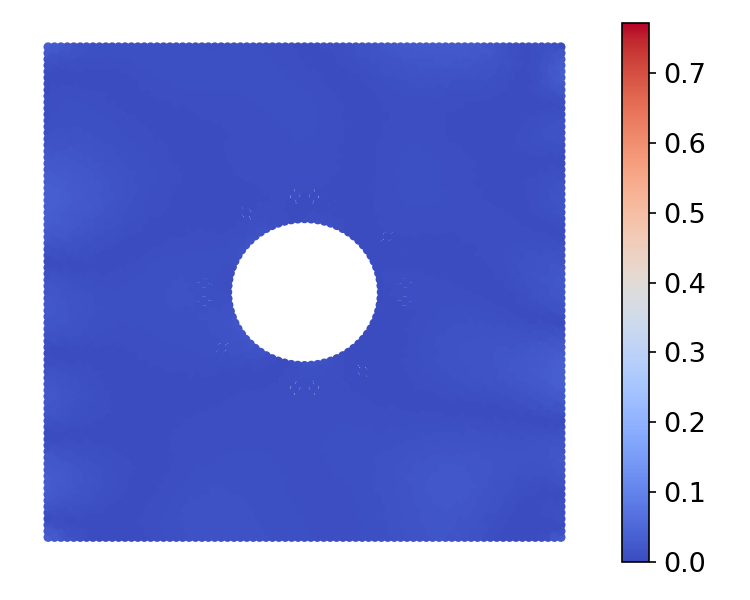}} \\ 
 & \rotatebox[origin=l]{90}{\makecell{PI-DeepONet \\ Worst case}}& {\includegraphics[width=0.16\textwidth]{ 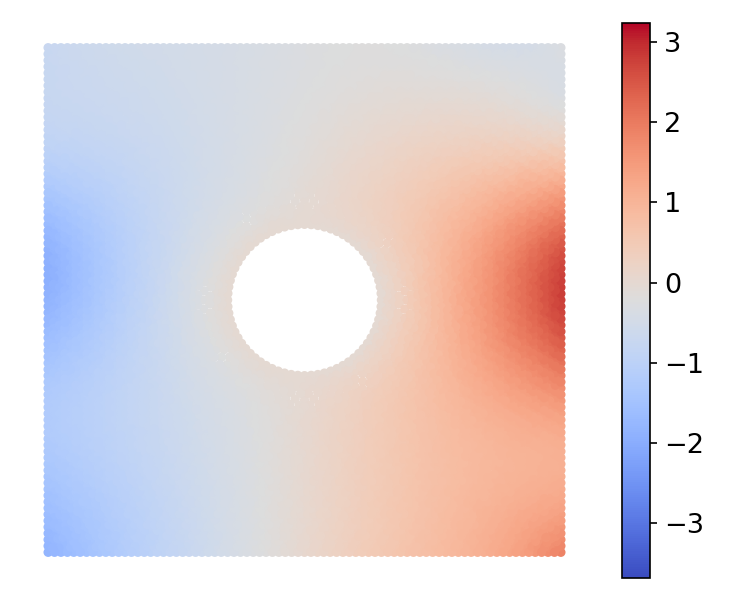}} & {\includegraphics[width=0.16\textwidth]{ 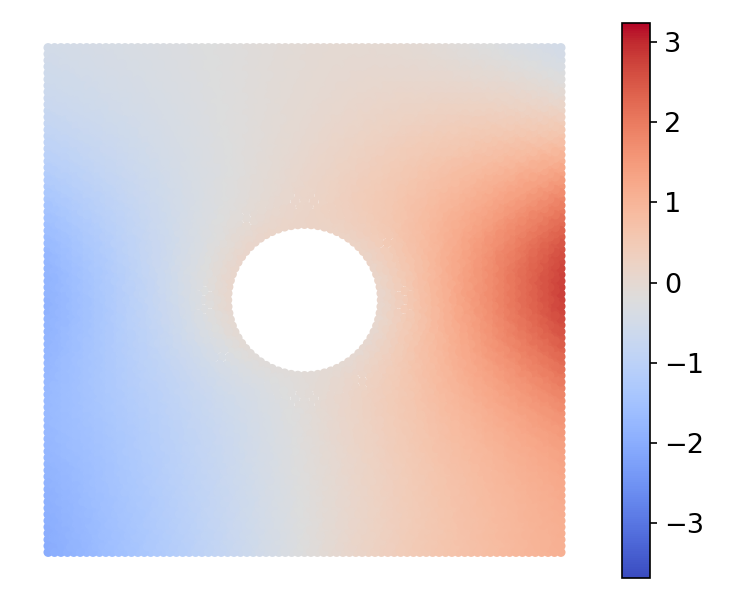}} & {\includegraphics[width=0.16\textwidth]{ 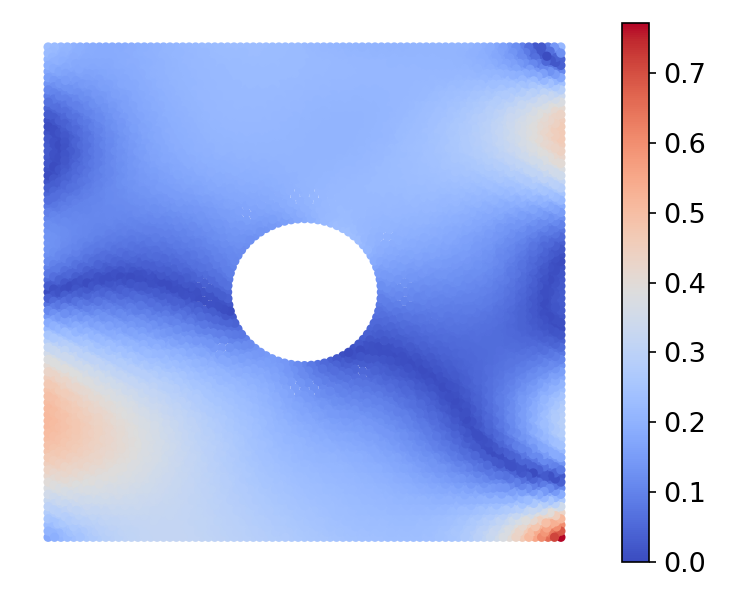}} & {\includegraphics[width=0.16\textwidth]{ 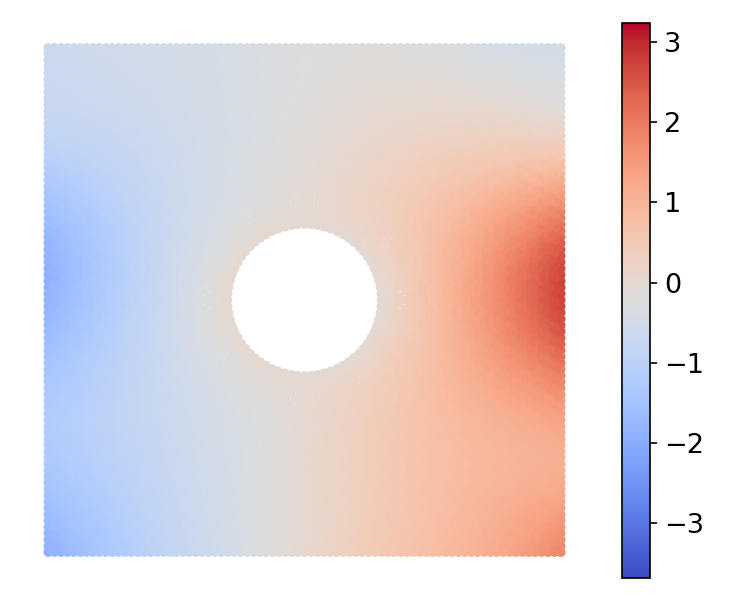}} & {\includegraphics[width=0.16\textwidth]{ 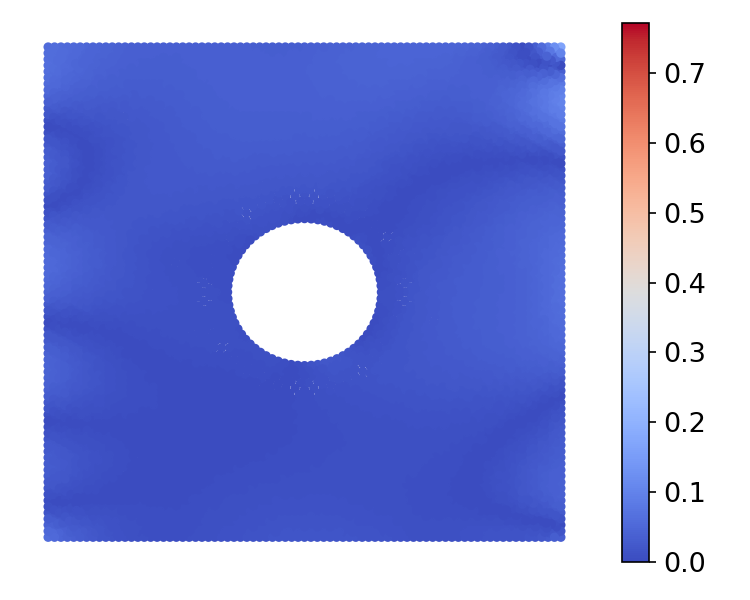}} \\
 & \rotatebox[origin=l]{90}{\makecell{PI-DCON \\ Worst case}} & {\includegraphics[width=0.16\textwidth]{ 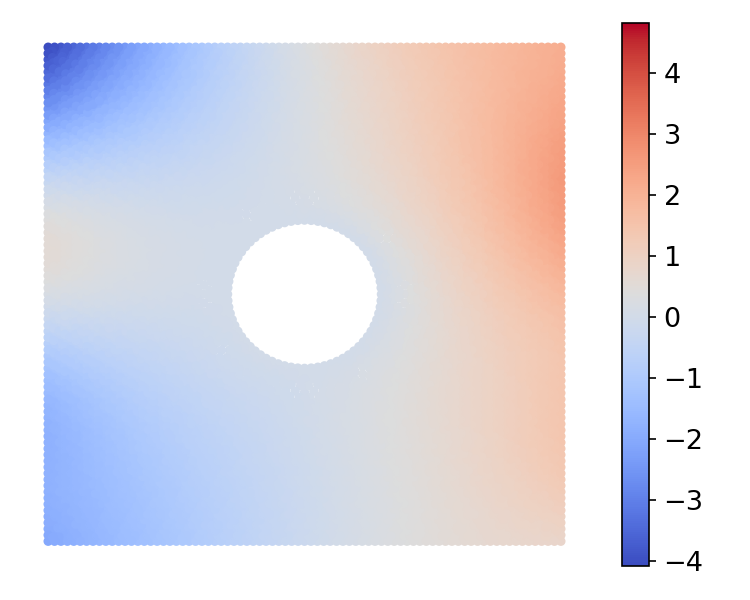}} & 
 {\includegraphics[width=0.16\textwidth]{ 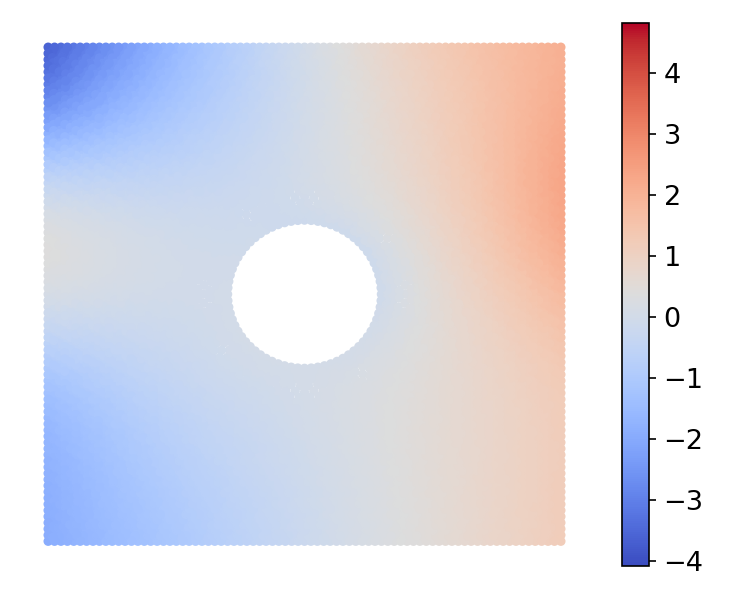}} & {\includegraphics[width=0.16\textwidth]{ 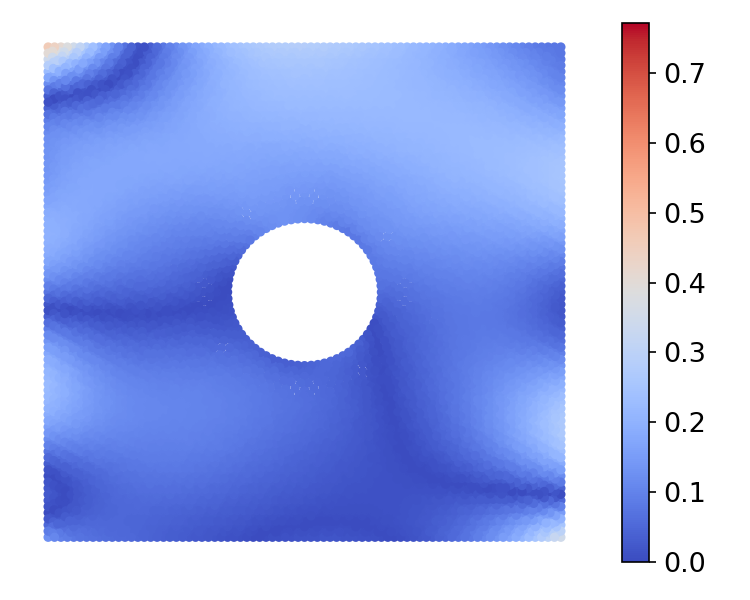}} &{\includegraphics[width=0.16\textwidth]{ 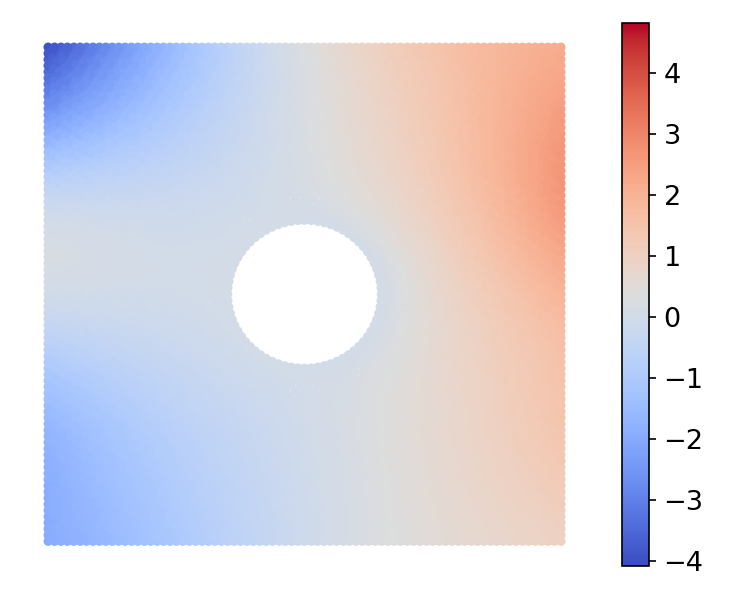}} & {\includegraphics[width=0.16\textwidth]{ 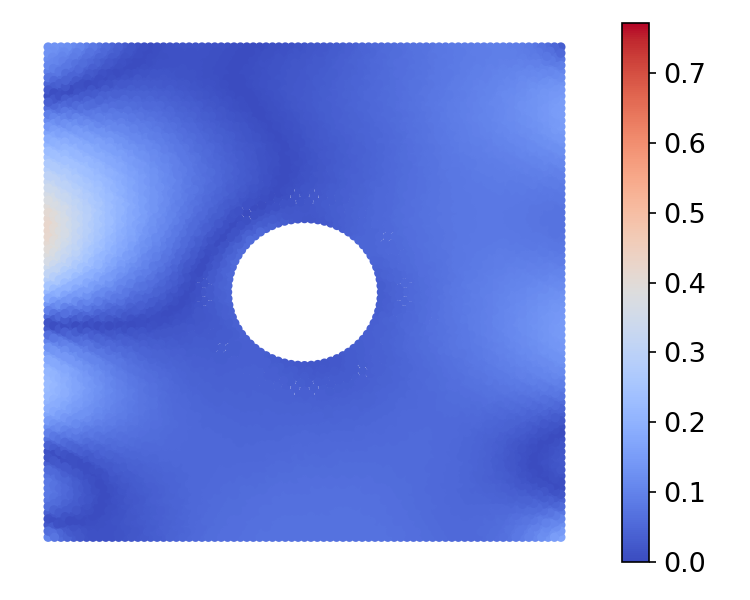}}  \\ 
 \hline
\end{tabular}
\label{fig.best_worst_case}
\end{table*}

\begin{table*}[ht]
\centering 
\caption{Performance of PI-DCON trained with variable discrete representations. The best and worst  prediction among various realizations of boundary conditions are shown.}
\begin{tabular}{|c|c|c c c|}
\hline
\multicolumn{2}{|c|}{} & Ground Truth & Prediction & Absolute Error \\ 
\hline
\multirow{2}{*}{\rotatebox[origin=c]{90}{Darcy flow}} & \rotatebox[origin=l]{90}{\makecell{Best case}} &  {\includegraphics[width=0.16\textwidth]{ 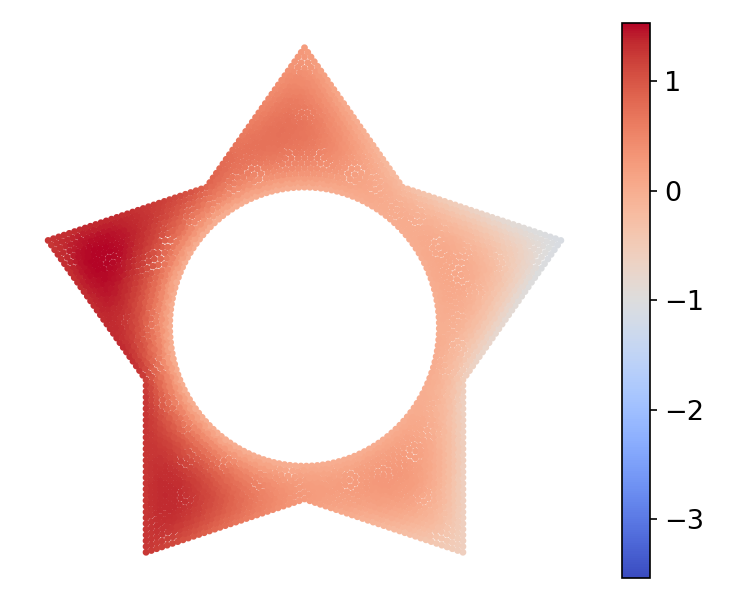}} & {\includegraphics[width=0.16\textwidth]{ 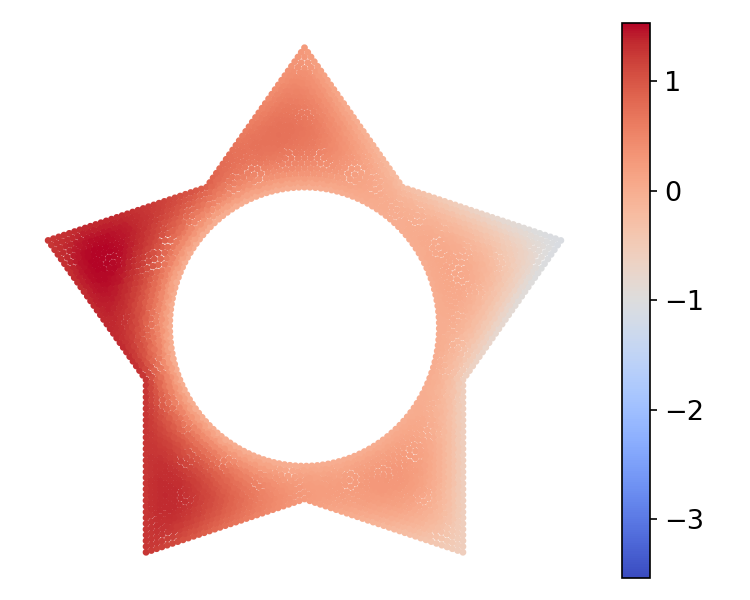}} & {\includegraphics[width=0.16\textwidth]{ 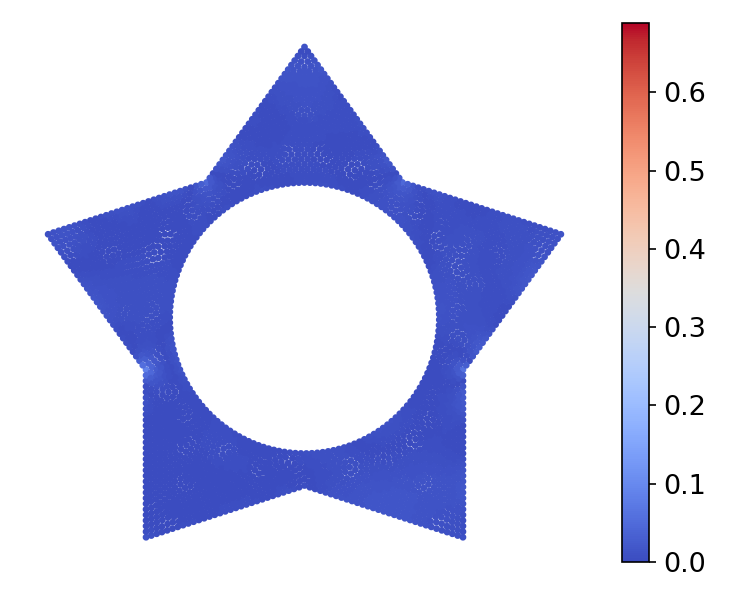}}\\ 
& \rotatebox[origin=l]{90}{\makecell{Worst case}} &  {\includegraphics[width=0.16\textwidth]{ 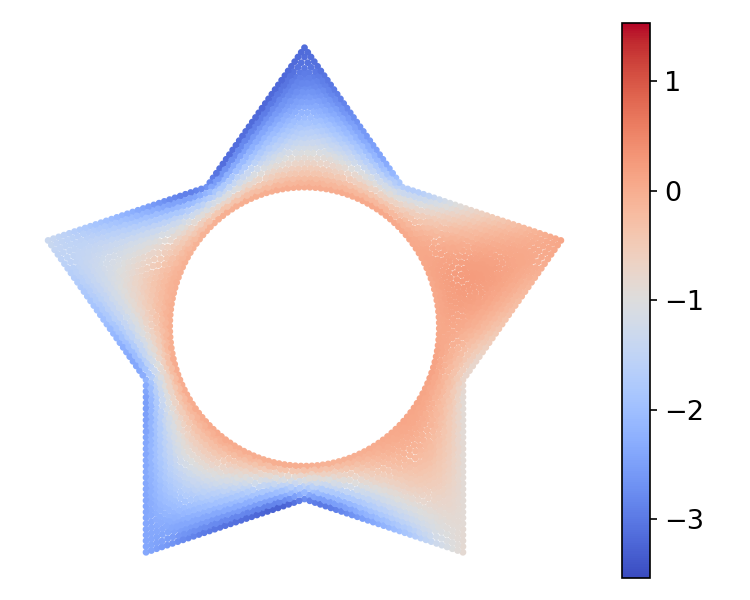}} & {\includegraphics[width=0.16\textwidth]{ 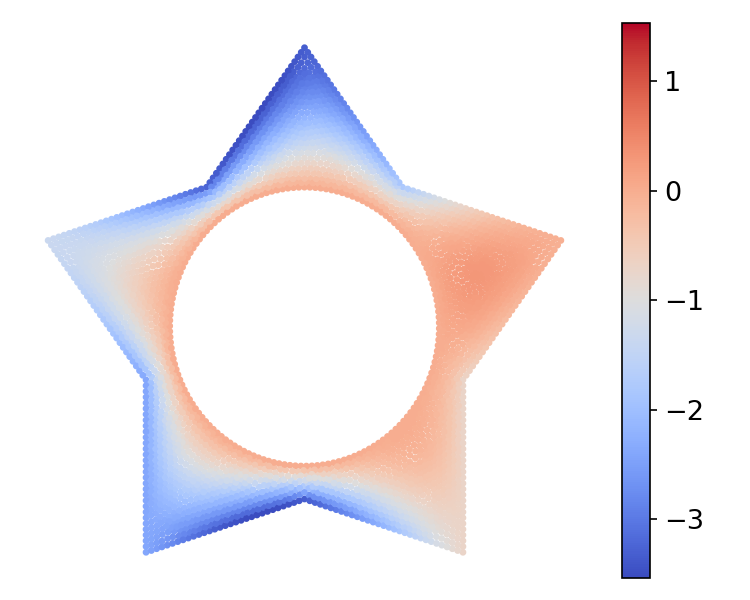}} & {\includegraphics[width=0.16\textwidth]{ 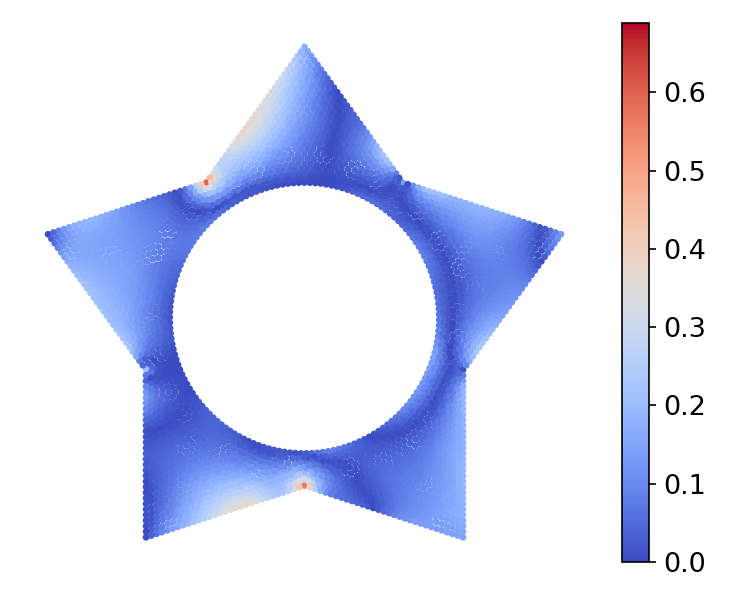}} \\ 
\hline
\multirow{2}{*}{\rotatebox[origin=c]{90}{2D plate}}
 & \rotatebox[origin=l]{90}{\makecell{Best case}} & {\includegraphics[width=0.16\textwidth]{ 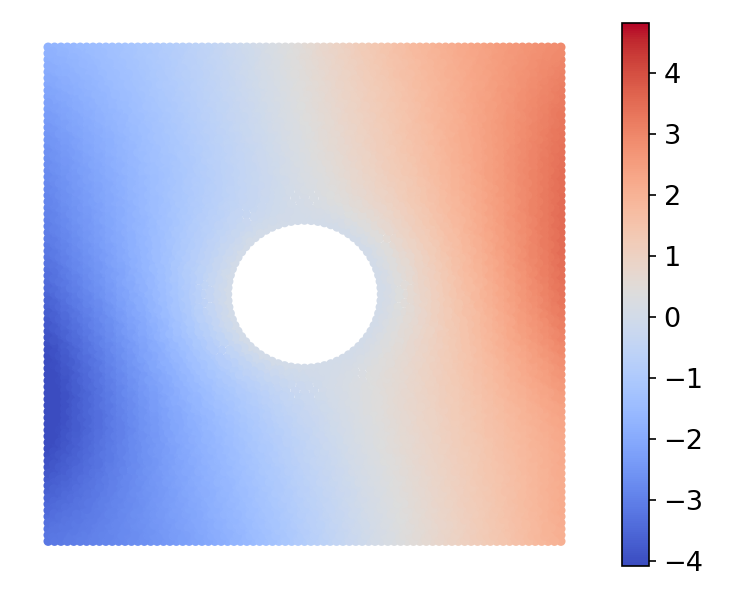}} & {\includegraphics[width=0.16\textwidth]{ 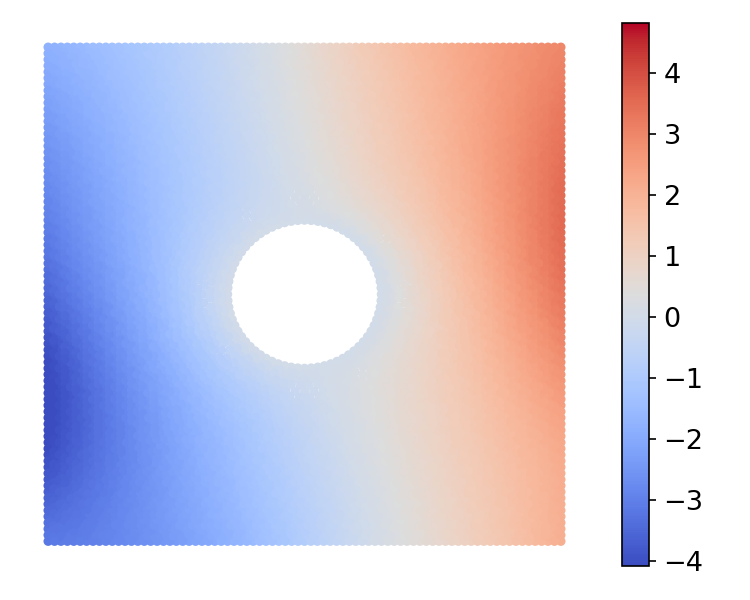}} & {\includegraphics[width=0.16\textwidth]{ 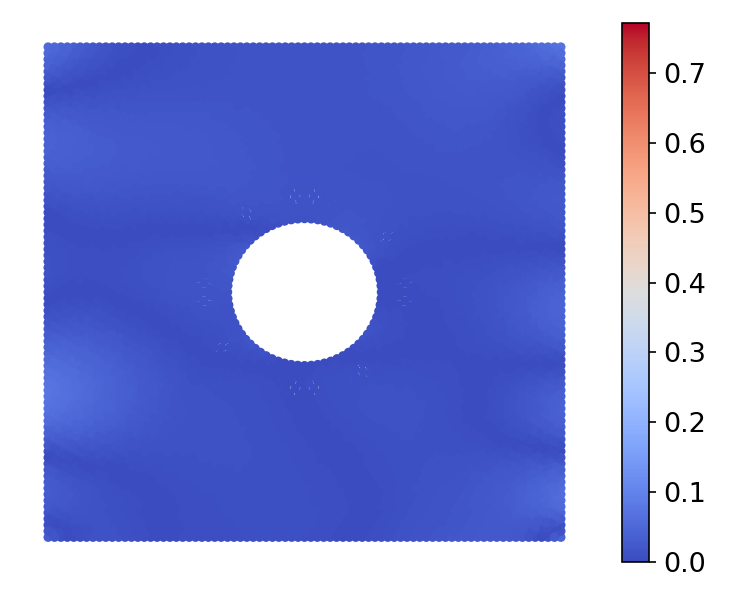}} \\
  & \rotatebox[origin=l]{90}{\makecell{Worst case}} & {\includegraphics[width=0.16\textwidth]{ 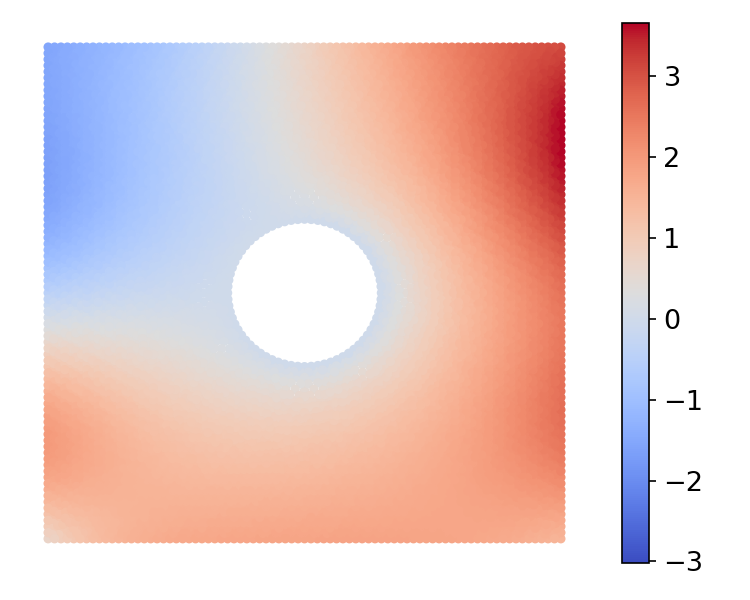}} & {\includegraphics[width=0.16\textwidth]{ 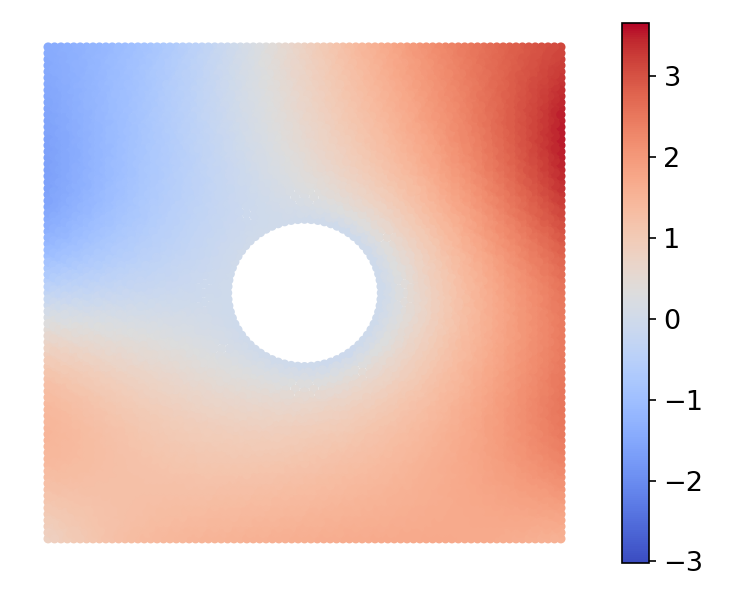}} & {\includegraphics[width=0.16\textwidth]{ 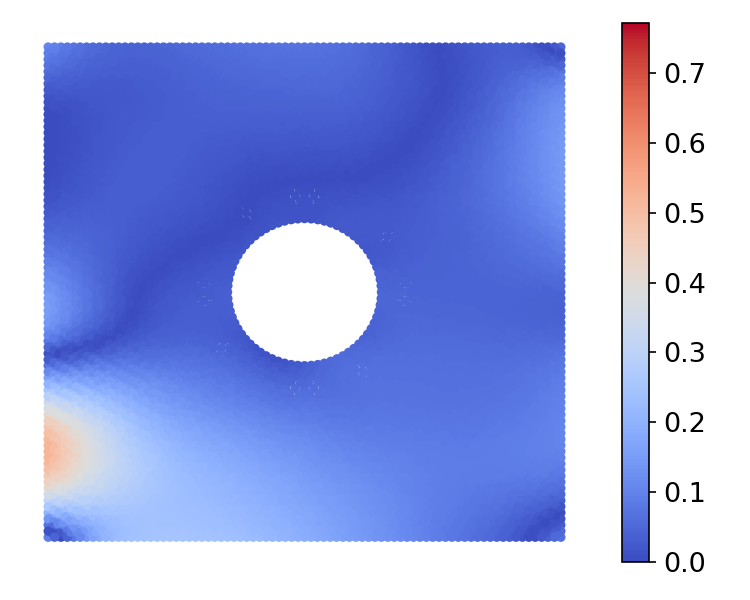}} \\
\hline
\end{tabular}
\label{fig.best_worst_case_vary_mesh}
\end{table*}

\subsection{Comparison with data-driven neural operators} \label{sec.datadriven}

One of the main features of PI-DCON is that it is a data-free approach. It doesn't require  FE simulation runs as a training set.  In this section, we seek to investigate its efficiency and accuracy in comparison with  data-driven (supervised) training approaches. In doing so, we consider the data-driven DeepONet and also create a data-driven variant of our proposed DCON architecture for comparison. We first compare the performance of data-driven DCON with that of data-driven DeepONet \textcolor{black}{data-driven version of \cite{IDON}} in Table \ref{table.data_driven_com}. \textcolor{black}{ We observe that the data driven DCON model achieved 3.8\% accuracy improvement for the Darcy flow problem and 20.7\% accuracy improvement for the 2D plate problem, compared with the best baseline method.} Also,  when we compare our model's performance on data with variable discrete representations to those with a fixed discrete representations, we see a small drop in the prediction accuracy, which is comparable to the drop observed in PI-DCON. 

\begin{table*}[!ht]
\begin{center}
\caption{Accuracy comparison between data-driven DCON and data-driven DeepONet. }
\begin{tabular}{c c c c}
\hline
\multirow{2}{*}{Dataset} & \multicolumn{3}{c}{\makecell{$L_2$ Relative error}} \\
& \makecell{Data-driven DeepONet} & \makecell{Data-driven \cite{IDON}} & \makecell{Data-driven DCON}\\
\hline
\makecell{Darcy flow \\ (fixed discrete representation)} & 3.74\% $\pm$ 0.97\% & 1.82\% $\pm$ 0.44\% & \bf{1.66\%} $\pm$ \bf{0.37\%}\\
\makecell{2D plate \\ (fixed discrete representation)} & 3.78\% $\pm$ 1.01\% & 1.98\% $\pm$ 0.41\% & \bf{1.57\%} $\pm$ \bf{0.31\%} \\
\hline
\makecell{Darcy flow \\ (variable discrete representation)} & - & - & 2.12\% $\pm$ 0.49\%\\
\makecell{2D plate \\ (variable discrete representation)} & - & - & 1.78\% $\pm$ 0.34\% \\
\hline
\label{table.data_driven_com}
\end{tabular}
\end{center}
\end{table*}

\begin{figure}[ht]
    \centering
            \begin{subfigure}[b]{0.3\textwidth}
            \centering
            \includegraphics[width=\textwidth]{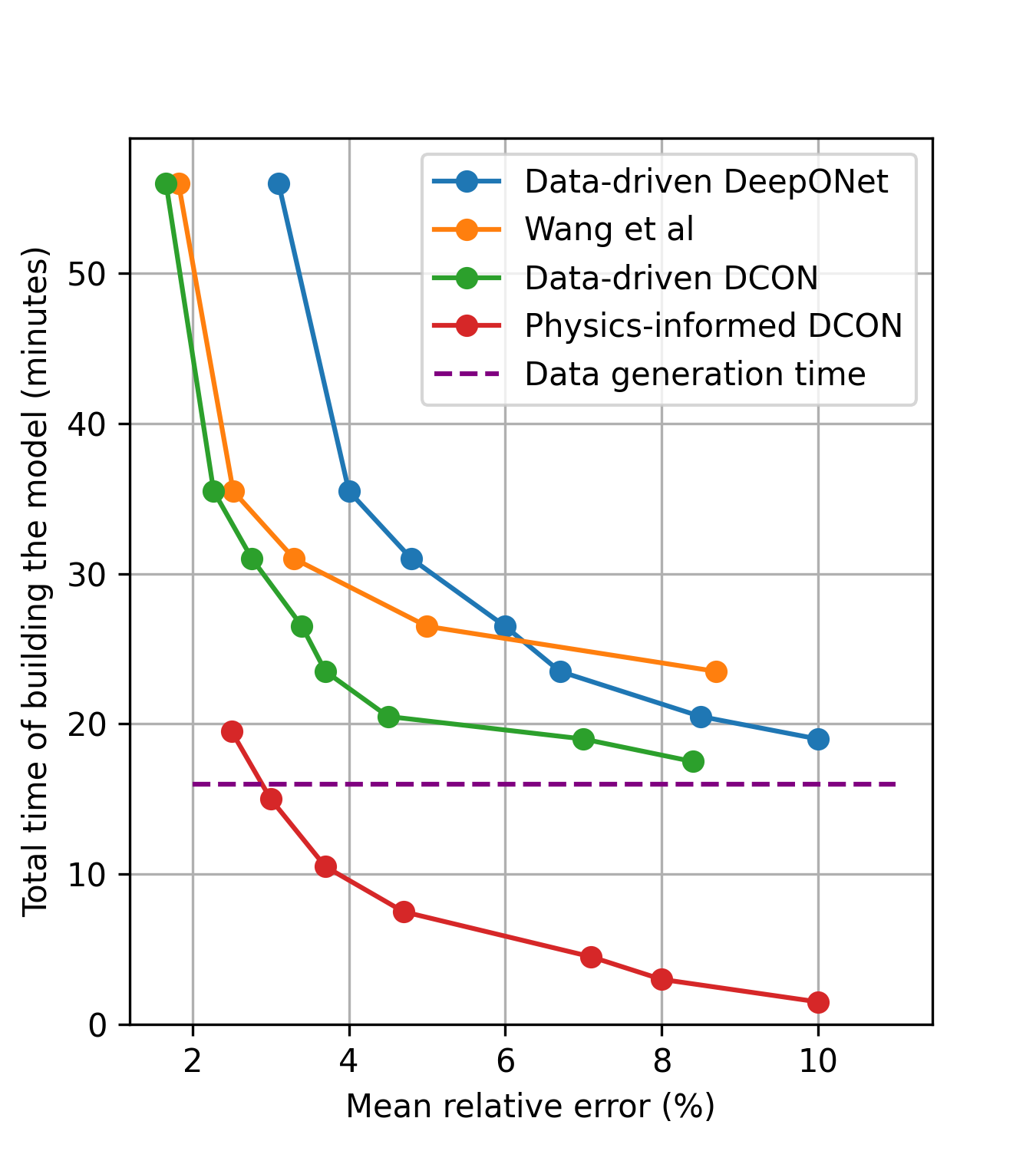}
            \caption{Darcy flow problem}    
            \label{fig:compare_DD_darcy}
        \end{subfigure}
        \hspace{0.1cm}
        \begin{subfigure}[b]{0.3\textwidth}  
            \centering 
            \includegraphics[width=\textwidth]{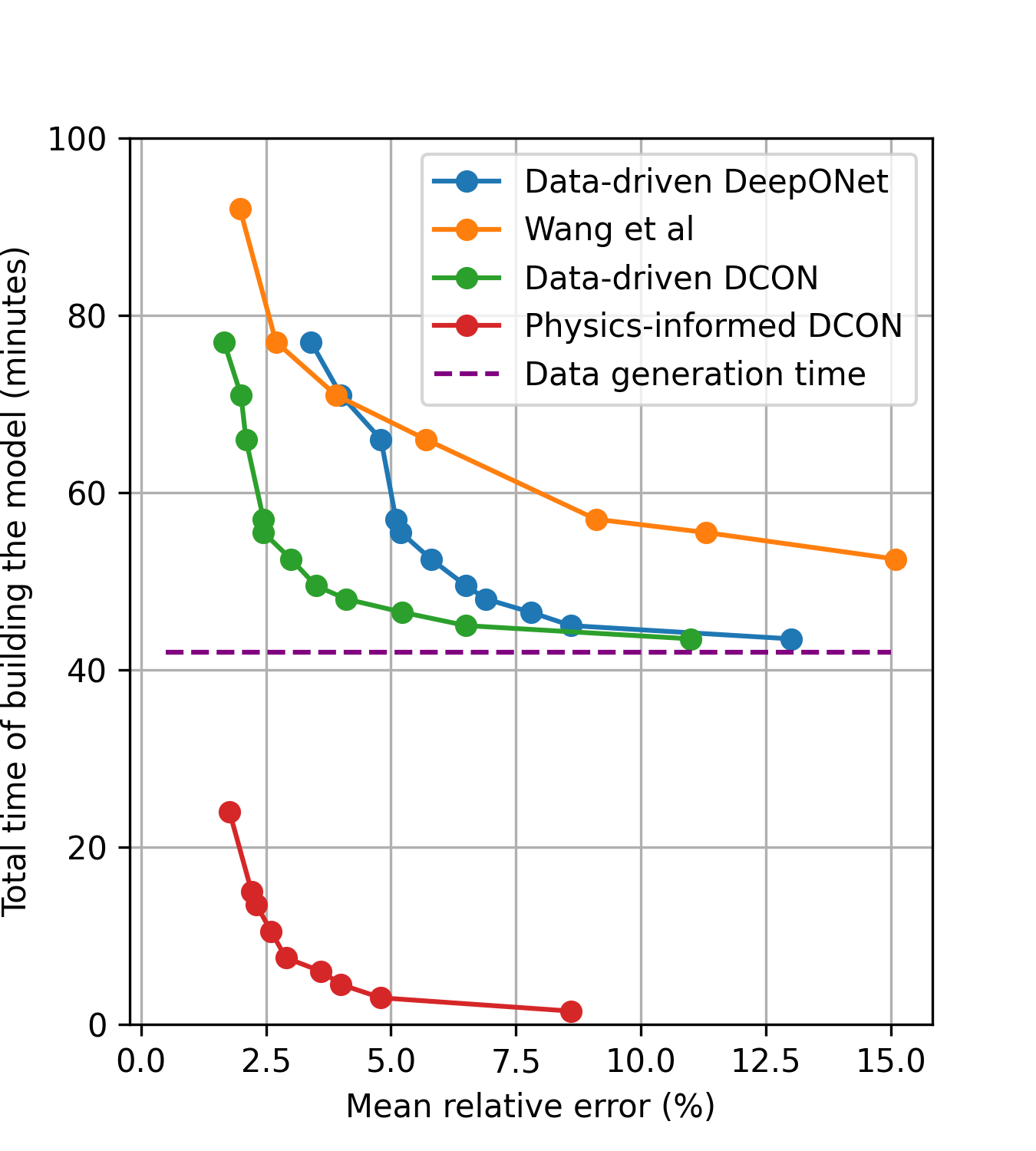}
            \caption{2D plate problem}    
            \label{fig:compare_DD_plate}
        \end{subfigure}
    \caption{\footnotesize \textcolor{black}{Comparison between the total time needed to train different neural operators. For data-driven method, this time also includes the time needed to generate the training (FE simulation) data.}}
    \label{fig.eff_acc_com}
\end{figure}

In order to demonstrate the advantage of physics-informed (data-free) training of DCON, in Figure \ref{fig.eff_acc_com} we show the required time needed for building  models at different accuracy levels. The horizontal (dashed) line  represents the required CPU time for generating training datasets for data-driven training of neural operators. In both darcy flow and 2d plate problems, the lowest error rate achieved by PI-DCON is smaller than that of data-driven DeepONet but greater than that of data-driven DCON. \textcolor{black}{Also, PIDCON can obtain comparable results with the improved DeepONet architecture of \cite{IDON}.} This  suggests that neural operators trained in a data-driven fashion tend to outperform those trained using physics-informed approaches when the same model architecture is used. However, it can be seen that PI-DCON convergence is superior compared to data-driven DeepONet \textcolor{black}{and data-driven version of \cite{IDON}}. This shows the enhanced architectural design in DCON  allows for faster convergence compared to the conventional DeepONet architectures. 

Moreover, PI-DCON demonstrates greater training efficiency. This is when we compare the  total time needed to train a model, including the time to obtain training data in the data-driven approaches. It was observed that PI-DCON required only 55\% of the total time needed by data-driven DCON in the Darcy flow problem and 31\% in the 2D plate problem. This efficiency gain suggests that adopting a physics-informed training approach can lead to significant computational time savings, especially for more complex problems. 

\subsection{Ablation studies \label{Subse.ablation}}

In this section, we conduct ablation studies on the proposed PI-DCON  model, with a particular emphasis on how the number of operator layers affects the  accuracy. To this end, we train PI-DCON models with varying numbers of operator layers (specifically,  1, 2, 3, and 4) for 500 epochs, and report the prediction accuracies in Table ~\ref{table.ablation_num_layer}. We observe that an increase in the number of operator layers leads to improved accuracy. However, the marginal improvement in prediction accuracy diminishes with each additional layer. For example, increasing the layer count from 1 to 2 yields an average accuracy improvement of 3.79\%, whereas the increment from 3 to 4 results in a mere 0.45\% improvement.

\begin{table*}[!ht]
\begin{center}
\caption{Mean of relative error of different model sizes.}
\begin{tabular}{c c c c c c c}
\hline
 & \multicolumn{6}{c}{Number of operator layers}\\
 \cline{2-7}
 & 1 & 2 & 3 & 4 & 5 & 6\\
\hline
\makecell{Darcy flow \\ (fixed discrete representation)} & 10.10\% & 4.20\% & 2.50 \% & 2.15\% & 2.06\%  & 2.03\% \\
\makecell{2D plate \\ (fixed discrete representation)} & 5.92\% & 2.57\% & 1.77\% & 1.56\% & 1.43\% & 1.38\% \\
\hline
\makecell{Darcy flow \\ (variable discrete representation)} & 12.58\% & 5.01\% & 3.42\% & 2.83\% & 2.66\% & 2.59\%\\
\makecell{2D plate \\ (variable discrete representation)} & 7.83\% & 4.52\% & 2.98\% & 2.34\% & 2.16\% & 2.07\% \\
\hline
\label{table.ablation_num_layer}
\end{tabular}
\end{center}
\end{table*}

\begin{figure}[ht]
    \centering
            \begin{subfigure}[b]{0.3\textwidth}
            \centering
            \includegraphics[width=\textwidth]{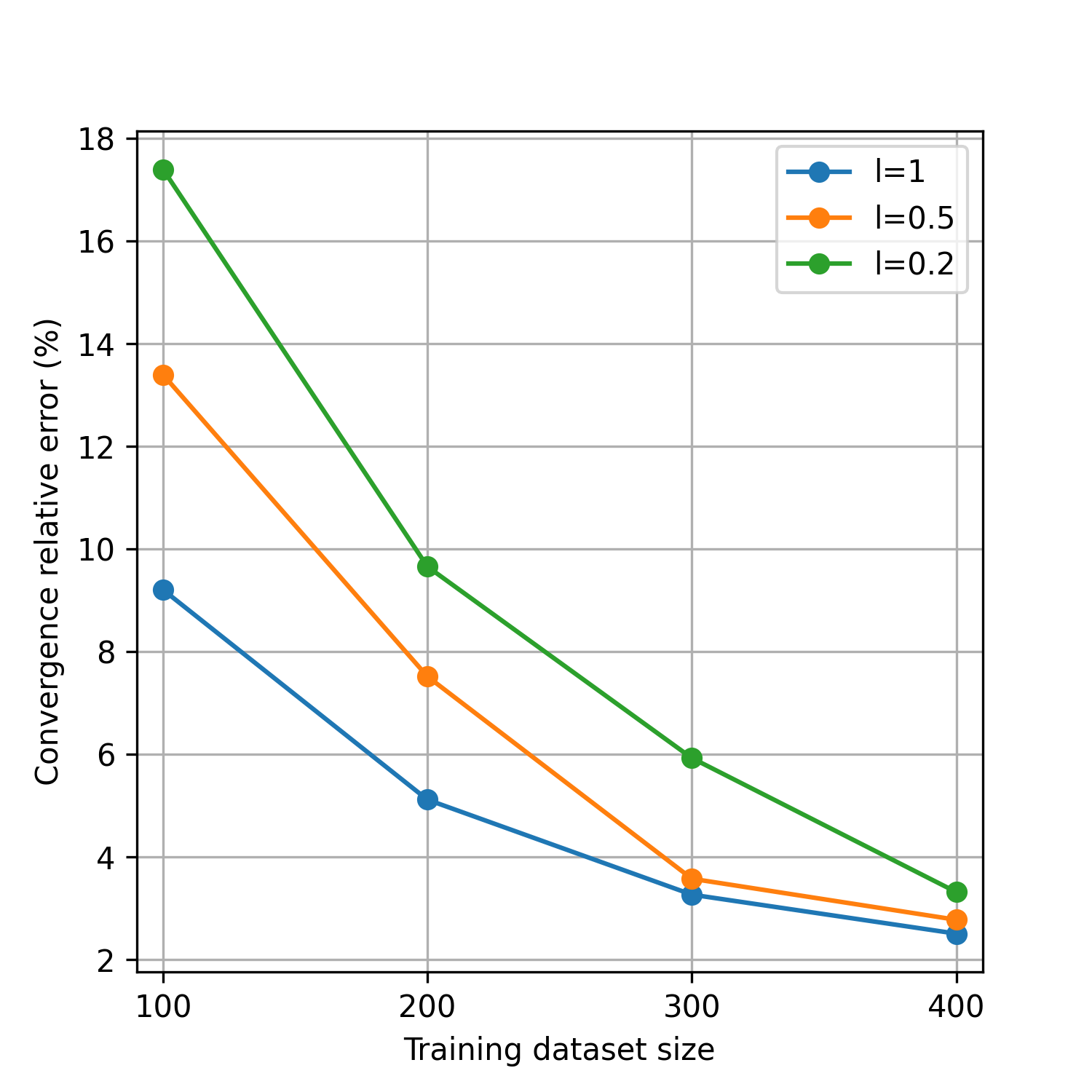}
            \caption{Darcy flow problem}    
            \label{fig:hyper_cl_darcy}
        \end{subfigure}
        \hspace{0.1cm}
        \begin{subfigure}[b]{0.3\textwidth}  
            \centering 
            \includegraphics[width=\textwidth]{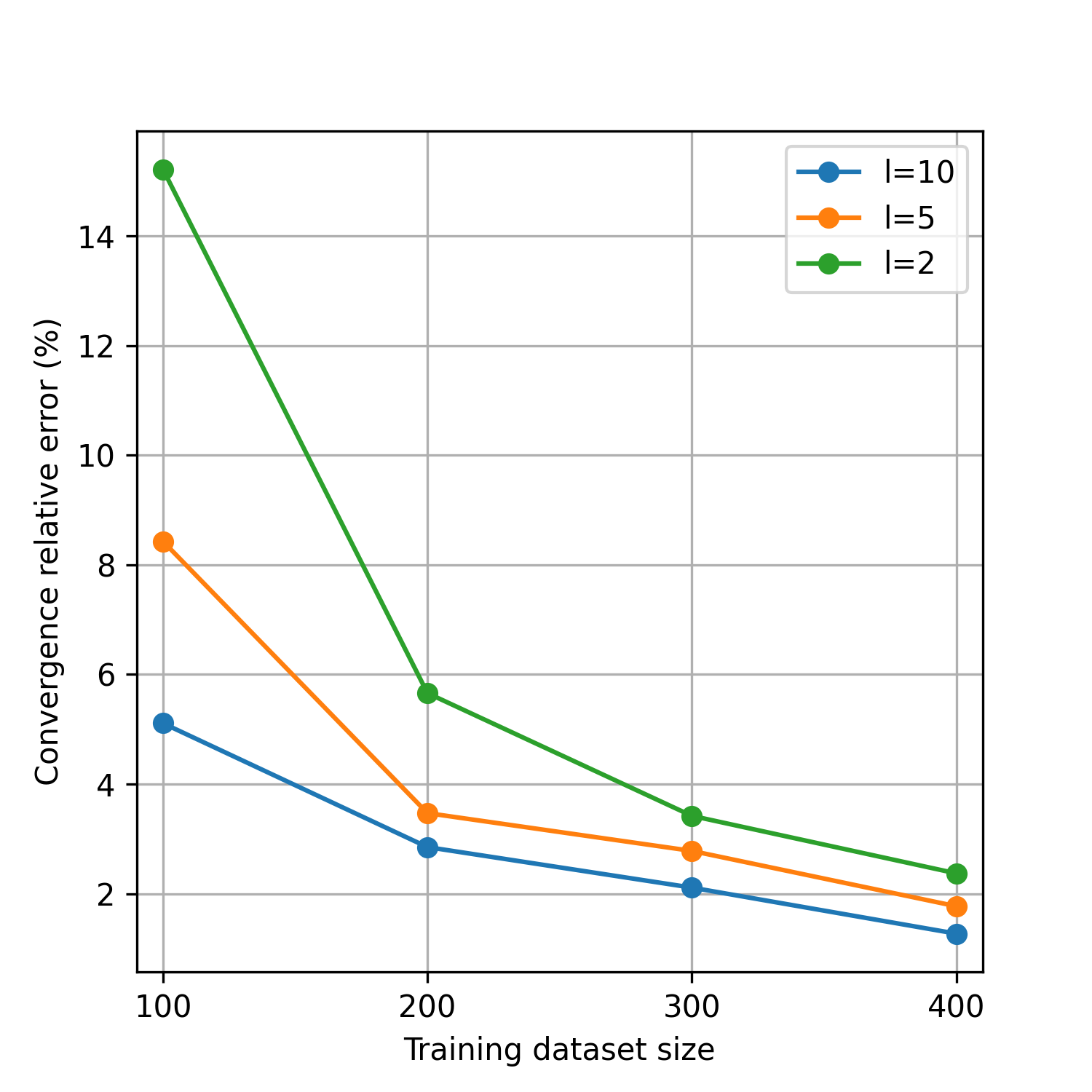}
            \caption{2D plate problem}    
            \label{fig:hyper_cl_plate}
        \end{subfigure}
    \caption{\footnotesize \textcolor{black}{Comparison between the total time needed to train different neural operators. For data-driven method, this time also includes the time needed to generate the training (FE simulation) data.}}
    \label{fig.ablation_cl}
\end{figure}

\textcolor{black}{We also seek to investigate our model performance with varying correlation lengths of the Gaussian Process and different amounts of training data. To do so, we train PI-DCON models using different numbers of training samples (specifically, 400, 300, 200, and 100 samples) for 500 epochs and report the prediction accuracy over 100 unseen samples with different correlation lengths $l$ in Figure ~\ref{fig.ablation_cl}. Our results show that increasing the number of training samples improves accuracy. For Gaussian processes with smaller correlation lengths, more training samples are required to achieve satisfactory performance. Additionally, if the training dataset size falls below a certain threshold, model performance significantly deteriorates. For instance, in the 2D plate problem, the relative error of the prediction increases from 5.66 to 15.22 when the dataset size decreases from 200 to 100, whereas it only increases from 3.42 to 5.66 when the dataset size decreases from 300 to 200.}

\section{Conclusion}\label{sec.conclusions}

In this study, we presented an enhanced architecture for physics-informed training of neural operators.  We introduced  PI-DCON based on the inspiration from  the Universal Approximation Theorem underlying the development of DeepONet. Our results indicate that PI-DCON achieves superior accuracy and generalization capabilities compared to PI-DeepONet. Furthermore, the comparison between our model with data-driven neural operators  highlights the advantages of developing neural operators using a physics-informed training approach.

Despite these advancements, our approach has its limitations. Primarily, the current model architecture is designed to generalize across various \textcolor{black}{ discrete representations} but does not extend this generalization to differing domain shapes. Furthermore, our evaluation focuses solely on the model's capability in computing steady-state solutions. Further research  is needed to extend PI-DCON to also incorporate dynamic responses, thereby enabling it to tackle time-dependent PDEs. \textcolor{black}{Also, the evaluation on other well-studied PDE problems, such as Eikonal equation and Stokes equation, can also be investigated.} \textcolor{black}{Additionally, while our model currently relies on the basic   multi-layer perceptron structure, incorporating more sophisticated architectures, such as Attention Mechanisms \cite{attn}, may further enhance its performance.} 

\textcolor{black}{In addition to developing more innovative architectures, investigating new physics-informed training algorithms can also provide options for improvement. As pointed out by \cite{IDON}, the balance between the PDE residual loss and the boundary condition loss significantly impacts model performance. Therefore, future works can focus on investigating various training algorithms that have shown promising results in PINN training, such as Variational PINNs (VPINNs) \cite{VPINN} and boundary distance function approximation \cite{ADF}.}

\bibliographystyle{plainnat}  
\bibliography{references}  

\end{document}